\documentclass[10pt]{amsart}
\RequirePackage[OT1]{fontenc}
\RequirePackage{amsthm,amsmath}
\RequirePackage[numbers]{natbib}
\RequirePackage[colorlinks,citecolor=blue,urlcolor=blue]{hyperref}
\usepackage{hhline}
\usepackage[table,xcdraw]{xcolor}
\usepackage{multirow}
\usepackage[toc,page,titletoc]{appendix}
\usepackage{xr-hyper}
\usepackage{hyperref}
\usepackage[british]{babel}
\usepackage{nameref}
\usepackage{amsfonts}
\usepackage{amssymb}
\usepackage{amsfonts}
\usepackage{mathrsfs}
\usepackage{amsfonts}
\usepackage{amstext}
\usepackage{amsopn}
\usepackage{pgfplots,geometry}
\usepackage{color}
\usepackage[]{subfig}
\usepackage[justification=justified]{caption}
\usepackage{tikz}
\usepackage{mhchem,gensymb}
\usepackage{graphicx}
\usepackage{chemfig}
\usepackage{url}

\usepackage{lmodern}
\usepackage{pifont}
\usepackage[abs]{overpic}

\usepackage{texdraw}
\usepackage{extarrows}
\usepackage{bbm}
\usepackage[utf8]{inputenc}
\usepackage{mathtools}

\usepackage{enumitem}

\allowdisplaybreaks
\usetikzlibrary{arrows, decorations.markings}
\usetikzlibrary{positioning}
\usetikzlibrary{fit}
\usetikzlibrary{backgrounds}
\usetikzlibrary{calc}
\usetikzlibrary{shapes}
\usetikzlibrary{mindmap}
\usetikzlibrary{decorations.text}
\pgfplotsset{compat=newest}

\usepgfplotslibrary{fillbetween}

\usetikzlibrary{intersections,quotes,angles}

\numberwithin{equation}{section}
\theoremstyle{plain}
\newtheorem{theorem}{Theorem}[section]
\newtheorem{lemma}[theorem]{Lemma}
\newtheorem{corollary}[theorem]{Corollary}
\newtheorem{proposition}[theorem]{Proposition}

\newtheorem{question}[theorem]{Question}

\theoremstyle{definition}

\newtheorem{example}[theorem]{Example}
\theoremstyle{remark}
\newtheorem{remark}[theorem]{Remark}
\newcommand{\eb}{\begin{example}}
\newcommand{\ee}{\end{example}}
\newcommand{\eqb}{\begin{equation}}
\newcommand{\eqe}{\end{equation}}
\newcommand{\spb}{\begin{split}}
\newcommand{\spe}{\end{split}}
\newcommand{\cab}{\begin{cases}}
\newcommand{\cae}{\end{cases}}
\newcommand{\thmb}{\begin{theorem}}
\newcommand{\thme}{\end{theorem}}
\newcommand{\qb}{\begin{question}}
\newcommand{\qe}{\end{question}}
\newcommand{\prob}{\begin{proposition}}
\newcommand{\proe}{\end{proposition}}
\newcommand{\leb}{\begin{lemma}}
\newcommand{\lee}{\end{lemma}}
\newcommand{\cob}{\begin{corollary}}
\newcommand{\coe}{\end{corollary}}
\newcommand{\enb}{\begin{enumerate}}
\newcommand{\ene}{\end{enumerate}}
\newcommand{\cenb}{\begin{center}}
\newcommand{\cene}{\end{center}}

\newcommand{\prb}{\begin{proof}}
\newcommand{\pre}{\end{proof}}
\newcommand{\rb}{\begin{remark}}
\newcommand{\re}{\end{remark}}

\newcommand{\tS}{\text{S}}
\newcommand{\ka}{\kappa}

\newcommand{\lt}{\left}
\newcommand{\rt}{\right}

\newcommand{\om}{\omega}

\newcommand{\N}{\mathbb{N}}
\newcommand{\bP}{\mathbb{P}}

\newcommand{\R}{\mathbb{R}}
\newcommand{\Z}{\mathbb{Z}}

\newcommand{\cP}{\mathcal{P}}

\newcommand{\cT}{\Omega}

\newcommand{\cX}{\mathcal{Y}}

\newcommand{\rO}{\mathrm{O}}
\newcommand{\ro}{\mathrm{o}}

\newcommand{\supp}{{\rm supp\,}}

\newcommand\sigman{\sigma_1}
\newcommand{\rhon}{\sigma_2}

\newcommand{\srcsize}{\@setfontsize{\srcsize}{5pt}{5pt}}

\makeatletter
\DeclareRobustCommand\widecheck[1]{{\mathpalette\@widecheck{#1}}}
\def\@widecheck#1#2{%
    \setbox\z@\hbox{\m@th$#1#2$}%
    \setbox\tw@\hbox{\m@th$#1%
       \widehat{%
          \vrule\@width\z@\@height\ht\z@
          \vrule\@height\z@\@width\wd\z@}$}%
    \dp\tw@-\ht\z@
    \@tempdima\ht\z@ \advance\@tempdima2\ht\tw@ \divide\@tempdima\thr@@
    \setbox\tw@\hbox{%
       \raise\@tempdima\hbox{\scalebox{1}[-1]{\lower\@tempdima\box
\tw@}}}%
    {\ooalign{\box\tw@ \cr \box\z@}}}
\makeatother
\begin{document}
\title[The asymptotic  tails of limit distributions of CTMCs]
  {The asymptotic tails of limit distributions of continuous time Markov chains}
\author{Chuang Xu}

\address{
Department of Mathematics\\
University of Hawai'i at M\={a}noa, Honolulu\\
96822, HI, US.}
\email{chuangxu@hawaii.edu}

\author{Mads Christian Hansen}

\address{
Department of Mathematical Sciences,
University of Copenhagen, Copenhagen,
2100 Denmark.}
\email{madschansen@gmail.com}

\author{Carsten Wiuf}

\address{
Department of Mathematical Sciences,
University of Copenhagen, Copenhagen,
2100 Denmark.}
\email{wiuf@math.ku.dk}


\date{\today}

\noindent

\begin{abstract}
 This paper investigates tail asymptotics of  stationary distributions and quasi-stationary distributions (QSDs) of continuous-time Markov chains on    subsets of the non-negative integers. {Based on the so-called flux-balance equation, we establish identities for stationary measures and QSDs,  which we use to derive   tail asymptotics.} In particular, continuous-time Markov chains with asymptotic power law transition rates, tail asymptotics for  stationary distributions and QSDs are classified into three types using three easily computable parameters: (i) super-exponential distributions, 
(ii) exponential-tailed distributions, and (iii) sub-exponential distributions. 
{Our approach to establish tail asymptotics of stationary distributions is different from the classical semimartingale approach, and we do not impose ergodicity nor moment bound conditions. In particular, the results also hold for explosive Markov chains, for which multiple stationary distributions may exist.} Furthermore, our results on tail asymptotics of QSDs seem new. 
We apply our results to biochemical reaction networks,  
a general single-cell stochastic gene expression model, an extended class of branching processes, and stochastic population processes with bursty reproduction, none of which are birth-death processes. The approach together with the identities easily extends to discrete time Markov chains.
\end{abstract}


\maketitle

\section{Introduction}

Stochastic biological models based on continuous-time Markov chains (CTMCs)  are commonly used to  model complex cellular behaviours, gene expression and the evolution of DNA \cite{E79}. Noise are inherent extrinsic and intrinsic properties of such biological systems and cannot be ignored without compromising the conclusions and the accuracy of the models.

In many cases it is reasonable to expect well-behaved biological systems, thus also well-behaved stochastic models. Therefore, in these cases, it is natural to assume the modelling CTMC is  \emph{ergodic}, that is, there exists a \emph{unique} stationary distribution which describes the system in the long run. In other cases, for example for  population processes without immigration, the population eventually goes extinct almost surely, and thus the ergodic stationary distribution is trivially the Dirac delta measure at zero. In these cases, it makes sense to study the \emph{quasi-stationary distribution} (QSD), that is,  the long-time behaviour of the process before extinction (usually called the \emph{$Q$-process}) \cite{MV12}. Jointly, stationary distributions and QSDs are referred to as \emph{limit distributions} in the present paper.

{The stationary distribution (provided it exists) is generally difficult to state in explicit form, except in  few cases. If the  underlying stochastic process has a \emph{detailed balanced} structure, then the stationary distribution takes a product-form;  see \cite{W86,Ke11,GM09,ACK10,MN10,CW16}  for various network models,  \cite{MV12} for birth-death processes (BDPs), and  \cite{G91,HM19,E05,N84,R88} for generalisations of such processes. QSDs with explicit expressions appear in even rarer cases \cite{V91}.}

While an explicit expression might not be known in general, less suffice in many cases. For example, if an expression for the  \emph{tail} distribution is known, then the existence and relative sizes of moments might be assessed from the decay rate of the tail distribution. Additionally, relative recurrence times might be assessed for stationary distributions of CTMCs.  

{With this in mind, we 
establish results for the \emph{tail} behaviour of   stationary distributions and QSDs, provided such exist. In particular, we concentrate on CTMCs on the non-negative integers $\N_0$  with asymptotic power law transition rate functions (to be made precise). }   
{Our approach is based on    generic identities   derived from the flux-balance equation \cite{Ke11} for limit distributions and stationary measures (Theorems~\ref{th-18} and \ref{th-18b}).
The identity for stationary distributions/measures  might be seen as a difference equation, which has order one less than the difference equation obtained directly from the master equation. More interestingly, for any given state $x$, the left hand side (LHS)  of the identity consists of  terms evaluated in states $\ge x$ only, while the right hand side (RHS)  of the identity consists of terms evaluated in states $<x$ only, and all terms have   \emph{non-negative} coefficients.    For BDPs, the identities coincide  with the classical recursive expressions for stationary distributions  and QSDs. }  
 
Furthermore, the identities allow us to  study the tail behaviour of limit distributions, provided they exist, and to   characterise their forms  (Theorems~\ref{th-19} and \ref{th-19b}). Specifically, in Section~5, for CTMCs with  transition rate functions that are asymptotically power law, we show that there are only three regimes: Either the decay follows   (i) a Conley-Maxwell-Poisson distribution (light-tailed), (ii) a geometric distribution, or (iii) a  sub-exponential distribution. Similar trichotomy result appears in the literature in other contexts \cite{A98,MZ06}, but to our knowledge, only for stationary distributions.  Our approach is based on repeated use of the identities we establish, combined with certain combinatorial identities (e.g., Lemma~\ref{Sle-1}).

Importantly,  we successfully obtain QSD tail asymptotics. To the   best of our knowledge, \emph{no similar classification results on QSD tail asymptotics have been established}, despite the Lyapunov function approach has been used frequently  to establish ergodicity of QSDs \cite{CV17}. A superficial reason might be    the inherent difference between stationary distributions and QSDs: Stationary distributions are equilibria of the master equation while QSDs are \emph{not}. A deeper explanation may be that the drift of  CTMCs (or DTMCs) may not directly relate to the decay rate of  QSDs. For an absorbing CTMC, known conditions for exponential ergodicity of stationary distributions are \emph{not} even sufficient to establish  uniqueness of QSDs \cite{V91}.

This difference is also  reflected in our results, where an extra condition is required to establish QSD tail asymptotics (Theorems~\ref{th-19b} and \ref{th-19}). Our novel approach successfully addresses both tail asymptoics of stationary distributions and QSDs at the same time, based on the similarity of the  algebraic equations they satisfy (Theorems~\ref{th-18} and \ref{th-18b}). 

We apply our main results to biochemical reaction networks, a general single-cell stochastic gene expression model, an extended class of branching processes, and stochastic population processes with bursty reproduction, none of which are BDPs.

{The Lyapunov function approach   is widely taken as a standard method to prove  ergodicity  of Markov processes \cite{MT09}. Additionally, the approach has been used to obtain  tail asymptotics of stationary distributions, assuming  ergodicity \cite{MP95,AI99,BGT01,DKW13,DKW16,XHW20b}. 
 In contrast, our results do not require ergodicity nor any other finite moment condition. For DTMCs, ergodicity follows from the existence of a stationary distribution  on an irreducible set \cite{M63}, while in contrast, this is not true for CTMCs \cite{M63}. 
 In particular, for \emph{explosive}  Markov chains, potentially with more than one stationary distribution, the cited results fail, whereas our results also hold in this case.}

{The trichotomy pattern for the tail asymptotics which we observe, is not surprising as it has already been observed for BDPs \cite{T04}, as well as for  processes on continuous state spaces; for example, the Lindley process    \cite{A98}, and the exponential functional of a L\'{e}vy process drifting to infinity  \cite{MZ06}. The techniques applied in these papers do not seem  applicable in our setting.}

It is noteworthy that the identities we establish, might be used to calculate the limit distribution \emph{recursively} up to an error term  that depends only on a few generating terms ($\pi(0)$ in the case of BDPs) and the truncation point of the limit distribution. The error term is given by the tail distribution, and thus, the decay rate of the error term might be inferred from the present work. Approximation of limit distributions will be pursued in a subsequent paper. The main results of this paper together with the approach can be extended to DTMCs in a straightforward way.

\section{Preliminaries}

\subsection{Sets and functions.}
Denote the set of real numbers, positive real numbers, integers, positive integers, and non-negative integers by $\R$, $\R_+$, $\mathbb{Z}$, $\mathbb{N}$, and $\mathbb{N}_0$, respectively. For $m, n\in\N$, let $\R^{m\times n}$ denote the set of $m$ by $n$ matrices over $\R$. Further, for any set $B$, let $\#B$ denote its cardinality and $\mathbbm{1}_B$ the corresponding indicator function. For $b\in\R$, $A\subseteq\R$, let $bA=\{ba\colon a\in A\}$, and $A+b=\{a+b\colon a\in A\}$. {Given $A\subseteq\R$, let $\inf A$ and $\sup A$ denote the infimum and supremum of the set, respectively. By convention, $\sup A=-\infty$, $\inf A=+\infty$ if $A=\varnothing$, and $\sup A=+\infty$ if $A$ is unbounded from above, and $\inf A=-\infty$ if $A$ is unbounded from below.} {For $x, y\in\mathbb{N}_0$ with $x\ge y$, let $x^{\underline{y}}=\frac{x!}{(x-y)!}$.} 

Let $f$ and $g$ be non-negative functions on an unbounded set $A\subseteq\mathbb{N}_0$. We denote $f(x)\lesssim g(x)$ if there exists $C,N>0$ such that
$$f(x)\le C g(x),\quad \text{for all}\quad x\in A,\ x\ge N,$$
that is, $f(x)=\rO(g(x))$ since $f$ is non-negative. Here O refers to the standard big O notation.
The function $f$ is said to be
\emph{asymptotic power law} (APL) if there exists $r_1\in\R$ such that  $\lim_{x\to\infty}\frac{f(x)}{x^{r_1}}=a$
 exists and is finite. Hence $r_1=\lim_{x\to\infty}\frac{\log f(x)}{\log x}$. An APL function $f$ is called \emph{hierarchical} (APLH) on $A$ with $(r_1,r_2,r_3)$ if there further exists $r_2, r_3$ with $r_2+1\ge r_1>r_2>r_3\ge r_1-2$, and $a>0$, $b\in\R$, such that for all large $x\in A$,
\begin{equation}\label{asymptotic-expansion}
f(x)=ax^{r_1}+bx^{r_2}+\rO(x^{r_3}).
\end{equation}
The requirement $r_2+1\ge r_1$ and $r_3\ge r_1-2$ comes from the analysis in Sections~6-7, where asymptotic Taylor expansion of functions  involves the  powers of the  first few leading terms. Here $r_1$, $r_2$, and $r_3$ are called the first, second and third power of $f$, respectively. All rational functions, polynomials, and real analytic APL functions are APLH. Not all APL functions are APLH, e.g., $f(x)=(1+(\log(x+1))^{-1})x$ on $\N$.

{Given a function $f$ which is APLH, the first power in the expansion is uniquely determined, while the other two powers are not. In other words, given the asymptotic expansion  \eqref{asymptotic-expansion}, there may exist a family of APLH functions admiting the same asymptotic expansion.} Let $r_2^*$ and $r_3^*$ be the infima over all $(r_2,r_3)\in\R^2$ such that $f$ is APLH on $A$ with $(r_1,r_2,r_3)$. For convention, we always choose the \emph{minimal} powers $(r_1,r_2^*,r_3^*)$, whenever $f$ is   APLH with $(r_1,r_2^*,r_3^*)$.
 As an example,
 $f(x)=x^2+3x+4$ is APLH on $\N_0$ with $(2,r_2,r_3)$ for any $2>r_2>r_3\ge 1$ (in which case $b=0$) or $1=r_2>r_3\ge 0$ ($b=3$). In this case, $f$ is APLH on $\N_0$ with minimal powers $(r_1,r_2^*,r_3^*)=(2,1,0)$. In contrast, take $f(x)=x+x^{1/3}\log x$. Then, $f$ is APLH on $\N_0$ with $(1,r_2,r_3)$ for any $r_2>r_3>1/3$ ($b=0$). In this case, $r_1=1$ and $r_2^*=r_3^*=1/3$, but $f$ is not APLH on $\N$ with $(1,1/3,1/3)$.
For any real analytic APLH function $f$ on $\N_0$, $f$ is APLH on $\N_0$ with $(r_1,r_2^*,r_3^*)$, where $r_1=\lim_{x\to\infty}\frac{\log f(x)}{\log x}$, $r_2^*=r_1-1$ and $r_3^*=r_1-2$.

\subsection{Markov chains.}

Let $(Y_t\colon t\ge 0)$ (or $Y_t$ for short) be a 
CTMC with state space $\cX\subseteq\N_0$ and transition rate matrix $Q=(q(x,y))_{x,y\in\cX}$, in particular each entry is finite. Recall that a set $A\subseteq\cX$ is \emph{closed} if $q(x,y)=0$ for all $x\in A$ and $y\in\cX\setminus A$ \cite{N98}.  Assume $\partial\subsetneq\cX$ is a (possibly empty)  finite closed \emph{absorbing set}, which is to the \emph{left} of $\partial^{\sf c}=\cX\setminus\partial$.  Here, the relative position of $\partial$ to $\partial^{\sf c}$ ensures that the only way for the Markov chain to end in an absorbing state is by jumping from a transitent state \emph{backward} to an absorbing state (this property is used in Proposition~\ref{prop-residual-eq} below).  
Furthermore, define {the set of jump vectors
$$\cT=\{y-x\colon q(x,y)>0,\ \text{for some}\ x,y\in\cX\},$$
and let $\cT_{\pm}=\{\omega\in\cT\colon {\rm sgn}(\om)=\pm1\}$ be the sets of forward and backward jump vectors, respectively.}  

For any probability distribution $\mu$ on $\partial^{\sf c}$, define $$\mathbb{P}_{\mu}(\cdot)=\int_{\cX}\mathbb{P}_x(\cdot){\rm d}\mu(x),$$
where $\mathbb{P}_x$ denotes the probability measure of $Y_t$ with initial condition $Y_0=x\in\cX$.
Any positive measure $\mu$ on a set $A\subseteq\N_0$ can be extended naturally  to  a positive measure on $\N_0$ with no mass outside $A$, $\mu(\N_0\setminus\! A)=0$. 

A (probability) measure $\pi$ on $\cX$ is a \emph{stationary measure} (distribution) of $Y_t$ if it is a non-negative equilibrium of the so-called {\em master equation} \cite{G92}:
 \eqb\label{Eq-24}0=\sum_{\omega\in\cT}q(x-\omega,x)\pi(x-\omega)-\sum_{\omega\in\cT}q(x,x+\omega)\pi(x),\quad x\in\cX.\eqe
(Here and elsewhere, functions defined on $\cX$ are put  to zero when evaluated at $x\not\in\cX\subseteq\Z$.)
Any stationary distribution $\pi$ of $Y_t$ satisfies 
$$\bP_{\pi}(Y_t\in A)=\pi(A),\quad A\in2^{\cX},$$
for all $t\ge0$, where $2^{\cX}$ is the power set of $\cX$ \cite{CV17}.

Let $\tau_{\partial}=\inf\{t>0\colon Y_t\in\partial\}$ be the entrance time of $Y_t$ into the absorbing set $\partial$. We say $Y_t$ admits \emph{certain absorption} if $\tau_{\partial}<\infty$ almost surely (a.s.) for all $Y_0\in\partial^{\sf c}$. Moreover, the process associated with $Y_t$ conditioned to  be never absorbed is called a {\em $Q$-process} \cite{CV16}.
A probability measure $\nu$ on $\partial^{\sf c}$ is a quasi-stationary distribution (QSD) of $Y_t$
 if for all $t\ge0$,
\[\mathbb{P}_{\nu}(Y_t\in A|\tau_{\partial}>t)=\nu(A),\quad  A \in 2^{\partial^{\sf c}}.\]

\section{Identities for limit distributions}

Let
$\omega_*=\text{gcd}(\Omega)$
be the (unique) positive greatest common divisor of $\cT$.
Define the scaled largest positive and negative jump, respectively:
\begin{align*}
\omega_+:=\sup_{\omega\in \Omega_+}\omega\om_*^{-1},\quad \omega_-:=\inf_{\omega\in \Omega_-}\omega\om_*^{-1}.
\end{align*}
Furthermore, define for $j\in\{\om_-,\ldots,\om_++1\}$,
\begin{align}\label{eq:Aj}
A_j &=\cab\{\omega\in\cT_-\colon j\omega_* >\omega\},\quad \text{if}\ j\in\{\omega_-,\ldots,0\},\\ \{\omega\in\cT_+\colon j\omega_*\le \omega\},\quad \text{if}\ j\in\{1,\ldots,\om_++1\}.\cae
\end{align}
 Hence, $\varnothing=A_{\om_-}\subseteq A_j\subseteq A_{j+1}\subseteq A_0=\cT_-$ for $\om_-< j<0$, and  $\varnothing=A_{\om_++1}\subseteq A_{j+1}\subseteq A_j\subseteq A_1=\cT_+$ for  $1< j< \om_+$.

The following classical result provides a necessary condition  for QSDs.

\prob[\rm \cite{CMM13}]\label{prop-residual-eq}
Assume $\partial\neq\varnothing$.
Let $\nu$ be a QSD of $Y_t$ on $\partial^{\sf c}$. Then for  $x\in\N_0\!\setminus\!\partial$,
\[
  \theta_{\nu}\nu(x)+\sum_{\omega\in\cT} q(x-\omega,x)\nu(x-\omega)-\sum_{\omega\in\cT}q(x,x+\omega)\nu(x) = 0,
\]
where
$$\theta_{\nu}=\sum_{\omega\in\cT_-}\sum_{y\in\partial^{\sf c}\cap\lt(\partial-\om\rt)}\nu(y)q(y,y+\omega)$$
 is finite.
\proe
\prb
 For $\partial=\{0\}$ and $x\in\partial^{\sf c}$, the identity follows from \cite{CMM13}. The same argument applies if $\partial^{\sf c}\not=\{0\}$.
For $x\in\N_0\setminus\cX$,  the identity trivially holds (with both sides being zero), which can be argued similarly to the proof of Theorem~\ref{th-18}.
\pre

 Before stating a generic identity for stationary distributions which originates from the master equation of the CTMC, we provide an example.

\begin{example}\label{WC}
Consider a CTMC on $\mathbb{N}_0$ with $\cT=\{1,-2\}$ and transition rate functions
\[q(x,x+1)=\kappa_1 x+\kappa_2 x(x-1),\quad q(x,x-2)=\kappa_3 x(x-1)(x-2),\]
where $\kappa_i$ ($i=1,2,3$) are positive constants. Such rate functions can be associated with a \emph{weakly reversible stochastic reaction network} (see Example~\ref{Ex-WC} in Section~\ref{sect-appl}). This CTMC is ergodic on $\mathbb{N}$ and there exists a stationary distribution $\pi$, which solves the master equation:
\begin{equation}\label{master-equation-identity}
b(x-1)\pi(x-1)+a(x+2)\pi(x+2)=(a(x)+b(x))\pi(x),\quad x\in\N_0,
\end{equation}
where $a(x)=\kappa_3x^{\underline{3}}$ and $b(x)=(\kappa_1x+\kappa_2x^{\underline{2}})$. One can show by induction that the terms in \eqref{master-equation-identity} might be separated such that   those with $\pi(y)$ for $y\ge x$ and those with $\pi(y)$ for $y<x$ are on different sides of the equality, in such a way that all coefficients are positive. Naive rearrangement of the terms does not suffice. We find
\begin{equation}\label{novel-identity}
a(x)\pi(x)+a(x+1)\pi(x+1)=b(x-1)\pi(x-1),\quad x\in\N_0,
\end{equation} 
 In addition, to our surprise, observe that \eqref{master-equation-identity} is a 3rd order difference equation, while   \eqref{novel-identity} is a 2nd order difference equation. Hence, the identity \eqref{novel-identity}, though equivalent to  \eqref{master-equation-identity}, helps simplify the relationship among all terms of the stationary distribution.
\end{example}

 The following generic identities generalise  the observations in Example~\ref{WC}, and  provide an equivalent definition  of a stationary distribution,    perhaps in a more handy form.  {The identities are the so-called flux-balance equation \cite[Lemma~1.4]{Ke11} specialized to one-dimension.} 
 We emphasise that stationary measures are not necessary finite, while stationary distributions and QSDs are probability distributions. 
\thmb \label{th-18}
The following statements are equivalent.
\begin{itemize}
\item[\rm{(1)}]  $\pi$ is a stationary measure (stationary distribution) of $Y_t$ on $\cX$.
\item[\rm{(2)}]  $\pi$ is a positive measure (probability distribution) on $\cX$ satisfying, for $x\in\N_0$,
\begin{equation*}
\sum_{\omega\in\cT_-}\sum_{j=\omega\omega^{-1}_*+1}^{0}q(x-j\omega_*, x-j\omega_*+\omega)\pi(x-j\omega_*)=
\sum_{\omega\in\cT_+}\!\sum_{j=1}^{\omega\omega^{-1}_*}
q(x-j\omega_*,x-j\omega_*+\omega)
\pi(x-j\omega_*)<\infty.
\end{equation*}
\item[\rm{(3)}] $\pi$ is a positive measure (probability distribution) on $\cX$ satisfying, for $x\in\N_0$,
\begin{equation*}
\sum_{j=\omega_-+1}^{0}\pi\lt(x-j\omega_*\rt)\sum_{\om\in A_j}q(x-j\omega_*,x-j\omega_*+\omega)
=  \sum_{j=1}^{\omega_+}\pi\lt(x-j\omega_*\rt)\sum_{\om\in A_j}q(x-j\omega_*,x-j\omega_*+\omega)<\infty.
\end{equation*}
\end{itemize}
\thme

\prb
We only prove for the equivalent representations for stationary measures. Then, the equivalent representations for stationary distributions follow. We assume w.l.o.g. that $\omega_*=1$.

{Recall from \cite[Lemma~1.4\footnote{{The identity \eqref{Stationary-Equation-A} is proved for finite $\mathcal{Y}$. It can be shown by induction that this identity holds true for countable $\mathcal{Y}$.}}]{Ke11} that $\pi$ is a stationary measure satisfying \eqref{Eq-24} if and only if for any $A\subseteq\mathcal{Y}$,
{\begin{equation}\label{Stationary-Equation-A}
\sum_{y\in A}\sum_{z\in A^c}\pi(y)q(y,z)=\sum_{y\in A}\sum_{z\in A^c}\pi(z)q(z,y)<\infty
\end{equation}
where $A^c=\mathcal{Y}\setminus A$. Letting $\pi(x)=0$ for $x\notin \mathcal{Y}$ and $q(x,y)=0$ for $(x,y)\notin\mathcal{Y}^2$, then \eqref{Stationary-Equation-A} holds for $A\subseteq\mathcal{Y}$ if and only if for $A\subseteq\Z$,
\begin{equation}\label{Stationary-Equation-A-extended}
\sum_{y\in A}\sum_{z\in \Z\setminus A}\pi(y)q(y,z)=\sum_{y\in A}\sum_{z\in\Z\setminus  A}\pi(z)q(z,y)<\infty
\end{equation}
In particular, for any $x\in \mathbb{N}_0$, let $A=\{y\in\Z\colon y\ge x\}$. Then \eqref{Stationary-Equation-A-extended} implies that
\begin{equation}\label{Stationary-Equation-x}
\sum_{y\ge x}\sum_{z\le x-1}\pi(y)q(y,z)=\sum_{y\ge x}\sum_{z\le x-1}\pi(z)q(z,y)<\infty
\end{equation}}
Conversely, \eqref{Stationary-Equation-x} implies \eqref{Eq-24} (i.e., $\pi$ is a stationary measure), which is shown by subtracting from both sides of
\eqref{Stationary-Equation-x} then same equations but with $x$ replaced by $x+1$. This is possible since both sides are  finite. We now compare  \eqref{Stationary-Equation-x} with Theorem~\ref{th-18}(3). {Recall that $A_{\om_-} = A_{\omega_++1} = \varnothing$.} In the latter equation, let $y=x-j$ and $z=x-j+\omega$, then \eqref{Stationary-Equation-x} follows from Theorem~\ref{th-18}(3). Furthermore, let $j=x-y$ and $\omega=z-y$, then \eqref{Stationary-Equation-x} is obtained from Theorem~\ref{th-18}(3). Observe that  Theorem~\ref{th-18}(2) and Theorem~ \ref{th-18}(3) are equivalent due to Fubini's theorem.
}
\pre

A special form of Theorem~\ref{th-18} under more assumptions has been stated in  the context of stochastic reaction   networks \cite[Prop.5.4.9]{H18}.



For any positive measure $\mu$ on $\N_0$, let
$$T_{\mu}\colon \N_0\to[0,1],\quad x\mapsto\sum_{y=x}^{\infty}\mu(y),$$
be the {\em tail distribution} (or simply the \emph{tail}) of $\mu$.

The following identities for QSDs also present equivalent definitions of the latter.

\thmb \label{th-18b}
Assume
$\partial\neq\varnothing$.
Then the following statements are equivalent.
\begin{itemize}
\item[\rm{(1)}] $\nu$ is a QSD of $Y_t$ on $\partial^{\sf c}$.
\item[\rm{(2)}] $\nu$ is a probability measure on $\partial^{\sf c}$, and for $x\in\N_0\setminus\partial$,
\begin{align*}&\sum_{\omega\in\cT_-}\sum_{j=\omega\omega_*^{-1}+1}^{0}q(x-j\omega_*,x-j\omega_*+\omega)\nu(x-j\omega_*)=\theta_{\nu}T_{\nu}(x)\\&+
\sum_{\omega\in\cT_+}\!\sum_{j=1}^{\omega\omega_*^{-1}}
q(x-j\omega_*,(x-j\omega_*+\omega)
\nu(x-j\omega_*)<\infty.\end{align*}
\item[\rm{(3)}]  $\nu$ is a probability measure on $\partial^{\sf c}$, and for $x\in\N_0\setminus\partial$,
\begin{align*}
&\sum_{j=\omega_-+1}^{0}\nu\lt(x-j\omega_*\rt)\sum_{\om\in A_j}q(x-j\omega_*,x-j\omega_*+\omega)
=\theta_{\nu}T_{\nu}(x)\\
&+\sum_{j=1}^{\omega_+}\nu\lt(x-j\omega_*\rt)\sum_{\om\in A_j}q(x-j\omega_*,x-j\omega_*+\omega)<\infty.
\end{align*}
\end{itemize}
\thme
\prb
The proof is similar to that of Theorem~\ref{th-18} and thus omitted.
\pre

{If a CTMC jumps uni-directionally (e.g., a pure birth or a pure death process), then all stationary measures, if such exist, are concentrated on {absorbing states} \cite{XHW20a}. In contrast, an absorbed pure {birth} process has no QSDs. However, an absorbed pure {death} process may admit one or more  QSDs, 
as illustrated below.}

\eb 
Consider a pure death process $Y_t$ on $\N_0$ with linear death rates $d_j=dj$ for $j\in\N_0$. Let $0<a\le 1$, and    
define $\nu$ as follows:
$$ \nu(1)=a,\quad \nu(x)=\left\{\begin{array}{cl}
 \frac{a}{x!}\frac{\Gamma(x-a )}{\Gamma(1-a)}, & \text{if}\ 0<a<1,\\ 
0, & \text{if}\ a=1,\end{array}\right\}\quad \text{for}\quad x>1, $$
where  $\Gamma(\cdot)$ is the Gamma function. Then, $\nu$ is  a QSD of $Y_t$ on $\N$ with $\partial=\{0\}$.    {Hence, there  is a family of QSDs fulfilling Theorem \ref{th-18b}.}
\ee

Formulae for stationary distributions and QSDs of BDPs follow directly from Theorems \ref{th-18} and \ref{th-18b}.

\cob[\cite{A91,C78}]\label{co-7}  
\rm{(i)} Let $Y_t$ be a BDP on $\N_0$ with birth and death rates $b_j$ and $d_j$, respectively, such that $b_{j-1}>0$ and $d_j>0$ for all $j\in\N$. If $\pi$ is a stationary distribution for $Y_t$, then \[\pi(j)=\pi(0)\prod_{i=0}^{j-1}\frac{b_i}{d_{i+1}},\quad j\in\N.\]

\rm{(ii)} Let $Y_t$ be a BDP on $\N_0$ with birth and death rates $b_j$ and $d_j$, respectively,  such that $b_0=0$, and $b_j>0$ and $d_j>0$ for all $j\in\N$. Then a probability distribution $\nu$ on $\N$ is a QSD trapped into $0$ for $Y_t$ if and only if
\eqb\label{Eq-2}d_j\nu(j)=b_{j-1}\nu(j-1)+d_1\nu(1)\lt(1-\sum_{i=1}^{j-1}\nu(i)\rt),\quad  j\ge2.\eqe
\coe

\prb
Here $\cT=\{-1,1\}$, $\om_*=\om_+=1$, $\om_-=-1$, and $\cX=\N_0$. Moreover, $q(j,j-1)=d_j$ and $q(j,j+1)=b_j$ for $j\in\N$.

\rm{(i)} $\partial=\varnothing$. Since $\pi$ is a stationary distribution on $\N_0$, it follows from Theorem~\ref{th-18} that
\[\pi(j)q(j,j-1)=\pi(j-1)q(j-1,j),\quad j\in\N.\]
 Hence the conclusion is obtained by induction.

\rm{(ii)}  $\partial=\{0\}$ and $\partial^{\sf c}=\N$. It follows from Theorem~\ref{th-18b} that a probability measure $\nu$ is a QSD on $\N$ if and only if
\[\theta_{\nu}=q(1,0)\nu(1),\quad \nu(j)q(j,j-1)=\theta_{\nu}T_{\nu}(j)+\nu(j-1)q(j-1,j),\quad j\in\N\setminus\!\{1\},\]
   that is, \eqref{Eq-2} holds.
\pre

Regarding the tail distributions, we have the following identities.

\cob\label{co-5}   Assume $\cT$ is finite and
$\partial=\varnothing$.
Let $\pi$ be a stationary distribution of $Y_t$  on $\cX$. Then, for $x\in\N_0$,
\begin{align}\nonumber
&T_{\pi}(x)\Bigl(\sum_{\omega\in A_0}q(x,x+\omega)+\sum_{ \omega\in A_1}q(x-\omega_*,x-\omega_*+\omega)\Bigr)+\sum_{j=\omega_-}^{-1}T_{\pi}(x-j\omega_*)\\ \nonumber
&\qquad\qquad \cdot\Bigl(\sum_{\om\in A_{j}}
q(x-j\omega_*,x-j\omega_*+\omega)-\sum_{\om\in A_{j+1}}
q(x-(j+1)\omega_*,x-(j+1)\omega_*+\omega)\Bigr)\\ \nonumber
&= \sum_{j=1}^{\omega_+}
T_{\pi}(x-j\omega_*)\Bigl(\sum_{\om\in A_j}
q(x-j\omega_*,x-j\omega_*+\omega)-\sum_{\om \in A_{j+1}}
q(x-(j+1)\omega_*,x-(j+1)\omega_*+\omega)\Bigr)
\end{align}
where $A_j$ is defined in \eqref{eq:Aj}.
\coe

\prb
Assume w.l.o.g.~$\omega_*=1$ and $0\in\cX$.
The LHS of the equation in Theorem~\ref{th-18}{\rm(3)} is
 \[\begin{split}
\text{LHS}=&\sum_{j=\omega_-+1}^{0}(T_{\pi}(x-j)-T_{\pi}(x-j+1))\sum_{\om\in A_j}  
q(x-j,x-j+\omega)\\
=&\sum_{j=\omega_-+1}^{0}T_{\pi}(x-j)\sum_{\om\in A_{j}}
q(x-j,x-j+\omega)\\
&-\sum_{j=\omega_-}^{-1}T_{\pi}(x-j)\sum_{\om\in A_{j+1}}
q(x-j-1,x-j-1+\omega)\\
=&\sum_{j=\omega_-}^{-1}T_{\pi}(x-j)\Big(\sum_{\om\in A_{j}}
q(x-j,x-j+\omega)-\sum_{\om\in A_{j+1}}
q(x-j-1,x-j-1+\omega)\Big)\\&+T_{\pi}(x)\sum_{\om\in A_0}
q(x,x+\omega),
\end{split}
\]while the RHS equals
\[\begin{split}
\text{RHS}=&\sum_{j=1}^{\omega_+}(T_{\pi}(x-j)-T_{\pi}(x-j+1))\sum_{\om\in A_j}
q(x-j,x-j+\omega)\\
=&\sum_{j=1}^{\omega_+}T_{\pi}(x-j)\Big(\sum_{\om\in A_j}
q(x-j,x-j+\omega)-\sum_{\om\in A_{j+1}}
q(x-j-1,x-j-1+\omega)\Big)\\&-T_{\pi}(x)\sum_{\om \in A_1}
q(x-1,x-1+\omega),
\end{split}
\]
which together yield the desired identity.
\pre

\cob
Assume $\cT$ is finite,
$\partial\neq\varnothing$, and let $\nu$ be a QSD of $Y_t$ on $\partial^{\sf c}$. Then for all $x\in\N_0\setminus\partial$,
\begin{align}\nonumber
&T_{\nu}(x)\Bigl(\sum_{\omega\in A_0}q(x,x+\omega)+\sum_{ \omega\in A_1}q(x-\omega_*,x-\omega_*+\omega)\Bigr)+\sum_{j=\omega_-}^{-1}T_{\nu}(x-j\omega_*)\\ \nonumber
&\qquad\qquad \cdot\Bigl(\sum_{\om\in A_{j}}
q(x-j\omega_*,x-j\omega_*+\omega)-\sum_{\om\in A_{j+1}}
q(x-(j+1),x-(j+1)+\omega)\omega_*)\Bigr)\\ \nonumber
&= \theta_{\nu}T_{\nu}(x)+\sum_{j=1}^{\omega_+}
T_{\nu}(x-j\omega_*)\Bigl(\sum_{\om\in A_j}
q(x-j\omega_*,x-j\omega_*+\omega)\\&-\sum_{\om \in A_{j+1}}
q(x-(j+1)\omega_*,x-(j+1)\omega_*+\omega)\Bigr),
\end{align}
where $A_j$ is defined in \eqref{eq:Aj}.
\coe

\begin{proof}
Similar to that of Corollary \ref{co-5}.
\end{proof}

\section{Asymptotic tails of limit distributions}\label{sect-Asymptotics}

To establish the asymptotic tails of limit distributions, we assume the following.

\medskip

\noindent($\rm\mathbf{A1}$) $\#\cT<\infty$.

\medskip

\noindent($\rm\mathbf{A2}$) $\cX$ is unbounded, and for  $\om\in\cT$, {$q(x,x+\omega)=a_{\omega}x^{R_{\om}^1}+b_{\omega}x^{R_{\om}^2}+\rO(x^{R_{\om}^3})$ is an APLH function in $x$ on $\cX$ with $(R_{\om}^1,R_{\om}^2,R_{\om}^3)$ for some constants $a_{\omega},\ b_{\omega}$. The APLH function is assumed} strictly positive for all large $x\in\cX$.

\medskip

\noindent($\rm\mathbf{A3}$) $\partial^{\sf c}$ is irreducible.

\medskip

Assumption ($\rm\mathbf{A1}$) guarantees that the chain   has  bounded jumps only, which  enable us to use the identities established in the previous section to estimate the tails.  
Assumption ($\rm\mathbf{A2}$) is common in applications, and ensures that  $\partial^{\sf c}$ is unbounded too, as $\partial$ is finite by assumption.
 In particular, ($\rm\mathbf{A2}$) is satisfied provided the following assumption holds:

\medskip

\noindent($\rm\mathbf{A2}$)' For $\om\in\cT$, $q(x,x+\omega)$ is a strictly positive polynomial for all large $x\in\cX$.

\medskip

Assumption  ($\rm\mathbf{A3}$) is assumed  to avoid non-essential technicalities. Moreover, ($\rm\mathbf{A3}$) means that either $Y_t$ is irreducible or the conditional process of $Y_t$ before entering $\partial$ is irreducible.
This assumption is satisfied for many known one-dimensional infinite CTMCs modeling biological processes (e.g., for population processes). In addition, ($\rm\mathbf{A3}$) implies that $\cT_+\neq\varnothing$ and $\cT_-\neq\varnothing$ (otherwise there are no non-singleton communicating classes).

The following parameters are well-defined and finite. Let
{$$R=\underset{{\om\in\cT}}{\max}R_{\om}^1,\quad R_-=\underset{{\om\in\cT_-}}{\max}R_{\om}^1,\quad R_+=\underset{{\om\in\cT_+}}{\max}R_{\om}^1,$$}
$$E_-=\underset{{\om\in\cT_-}}{\cup}\{R_{\om}^1, R_{\om}^2, R_{\om}^3\},\quad E_+=\underset{{\om\in\cT_+}}{\cup}\{R_{\om}^1, R_{\om}^2, R_{\om}^3\},$$
and define
$$\sigman=\min\{R_--R_-^1,R_+-R_+^1\},\quad \rhon=\min\{R_--R_-^2,R_+-R_+^2\},$$
where
$$R_-^1=\max \{r\in E_-\colon r<R_-\},\quad R_+^1=\max \{r\in E_+\colon r<R_+\},$$
 $$R_-^2=\max \{r\in E_-\colon r<R_-^1\},\quad R_+^2=\max \{r\in E_+\colon r<R_+^1\}.$$
(These values are only used to define $\sigman$, $\rhon$.)

Hence $R^1_-\ge R_{\om}^2\ge R_--1$, for some $\om\in\cT_-$ with {$R_{\om}^1=R_-$}. Similarly, $R^1_+\ge R_+-1$, $R^2_-\ge R_--2$, and $R^2_+\ge R_+-2$. This implies that $0<\sigman\le1$ and $\sigman<\rhon\le2$. If all transition rate functions are \emph{real analytic}, then by convention, $R_{\om}^2=R_{\om}^1-1,\ R_{\om}^3=R_{\om}^1-2$, for all $\om\in\cT$, and hence $\sigman=1$ and $\rhon=2$.
 Furthermore, let
 $$\alpha=\lim_{x\to\infty}\frac{\sum_{\omega\in\cT}q(x,x+\omega)\omega}{x^R},\quad \alpha_-=\lim_{x\to\infty}\frac{\sum_{\omega\in\cT_-}q(x,x+\omega)\vert \omega\vert}{x^{R_-}},$$
$$\alpha_+=\lim_{x\to\infty}\frac{\sum_{\omega\in\cT_+}q(x,x+\omega)\omega}{x^{R_+}},\quad \beta=\lim_{x\to\infty}\frac{\sum_{\omega\in\cT_-}q(x,x+\omega)}{x^{R_-}},$$
$$\gamma=\lim_{x\to\infty}\frac{\sum_{\omega\in\cT}q(x,x+\omega)\omega-\alpha x^R}{x^{R-\sigman}},\quad \vartheta=\frac{1}{2}\lim_{x\to\infty}\frac{\sum_{\omega\in\cT}q(x,x+\omega)\omega^2}{x^R}.$$
{The parameters also admit limit-free representations:
\begin{align*}
\alpha=\sum_{\omega\in\Omega\colon R_{\omega}^1=R}a_{\omega}\omega,\quad \alpha_-=\sum_{\omega\in\Omega_-\colon R_{\omega}^1=R_-}a_{\omega}\omega,\quad \alpha_+=\sum_{\omega\in\Omega_+\colon R_{\omega}^1=R_+}a_{\omega}|\omega|\\
\beta=\sum_{\omega\in\Omega\colon R_{\omega}^1=R}a_{\omega},\quad \gamma=\sum_{\omega\in\Omega\colon R_{\omega}^1=R-\sigma_1}a_{\omega}\omega+\sum_{\omega\in\Omega\colon R_{\omega}^2=R-\sigma_1}b_{\omega}\omega,\quad \vartheta=\sum_{\omega\in\Omega\colon R_{\omega}^1=R}a_{\omega}\omega^2.
\end{align*}
The  form with the limit  emphasizes that the parameters are coefficients of monomials of certain leading degrees. 
} Furthermore, define
 $$\Delta=\left\{\begin{array}{cl} -\gamma (\alpha_+\om_*)^{-1}, & \text{if}\quad \sigman<1,\\
 (-\gamma+R\vartheta) (\alpha_+\om_*)^{-1},& \text{if}\quad \sigman=1,\end{array}\right. \quad \delta=\Delta(\om_+-\om_--1)^{-1}.$$
{Note that $\alpha\le0$ implies $R_-\ge R_+$.} Moreover, $\alpha_+,\ \alpha_->0$ and $0<\delta\le\Delta$, due to ($\rm\mathbf{A3}$). Furthermore,  $\alpha, \alpha_-, \alpha_+, \beta$ and $\vartheta$ do not depend on the choice of second and third powers of the transition rate functions, whereas $\sigman$, $\rhon$, $\gamma$, $\Delta$ and $\delta$ do depend on the powers.

To state the results on the asymptotic tails of limit distributions, we classify probability distributions into the following classes.

 Let $\cP$ be the set of probability distributions on $A$. For  $a,b>0$, define
\begin{align*}
\cP_{a}^{1+}&=\{\mu\in\cP\colon T_{\mu}(x)\lesssim \exp(-ax(\log x) (1+\ro(1)))\},\\
 \cP_{a}^{1-}&=\{\mu\in\cP\colon T_{\mu}(x)\gtrsim \exp(-ax(\log x) (1+\ro(1)))\},\\
 \cP_{a,b}^{2+}&=\{\mu\in\cP\colon T_{\mu}(x)\lesssim \exp(-bx^a(1+\ro(1)))\},\\
  \cP_{a,b}^{2-}&=\{\mu\in\cP\colon T_{\mu}(x)\gtrsim \exp(-bx^a(1+\ro(1)))\},\\
\cP_{a}^{3+}&=\{\mu\in\cP\colon T_{\mu}(x)\lesssim x^{-a}\},\\ \cP_{a}^{3-}&=\{\mu\in\cP\colon T_{\mu}(x)\gtrsim x^{-a}\},
\end{align*}
where { $T_{\mu}(x)$ is the tail distribution of  a probability measure $\mu$,} and  $\ro(\cdot)$ refers to  the standard little o notation. Furthermore, define
\begin{align*}
\cP_{a}^{2+}=\cup_{b>0}\cP_{a,b}^{2+},&\hspace{1cm}  \cP_{a}^{2-}=\cup_{b>0}\cP_{a,b}^{2-},&\hspace{-4cm}  \cP_{<1}^{2-}=\cup_{0<a<1}\cP_{a}^{2-},\\
 \cP^{i+}=\cup_{a>0}\cP_{a}^{i+},  &\hspace{1cm} \cP^{i-}=\cup_{a>0}\cP_{a}^{i-}, \quad i=1,2,3.
\end{align*}
The sets  $\cP_{a}^{i+}$,  $i=1,2,3$,   are decreasing in $a$, while $\cP_{a}^{i-}$,  $i=1,2,3$, are increasing in $a$. Similarly, $\cP_{a,b}^{2+}$  is decreasing in both $a$ and $b$, while $\cP_{a,b}^{2-}$ is increasing in both $a$ and $b$. Morevoer, it is readily verified that
\[ \cP^{1+}\subseteq\cP_1^{2+}\subseteq\cP^{3+},\quad \cP^{3-}\subseteq\cP^{2-}_{<1}\subseteq\cP^{2-}_1\subseteq\cP^{1-}. \]
The probability distributions in $\cP^{2+}_{1}\cap\cP^{2-}_{1}$ decay as fast as exponential distributions and are therefore \emph{exponential-like}.  
Similarly, those in $\cP^{1+}$ are \emph{super-exponential};   
those in $\cP^{2-}_{<1}$ are \emph{sub-exponential}, and in particular those in  $\cP^{3+}\cap\cP^{3-}$ are \emph{power-like}.  
\cite{JKK05}.

The  \emph{Conley-Maxwell-Poisson} (CMP) distribution on $\N_0$ with parameter  $(a,b)\in\R_+^2$ has   probability mass function  given by \cite{JKK05}:
\[{\sf CMP}_{(a,b)}(x)=\frac{a^x}{(x!)^b}\lt(\sum_{j=0}^{\infty}\frac{a^j}{(j!)^b}\rt)^{\!\!\!-1},\quad x\in\N_0.\]
In particular, ${\sf CMP}_{a,1}$ is a Poisson distribution. For every probability distribution $\mu\in\cP^{1+}\cap\cP^{1-}$, there exists $(a_1,b_1),  (a_2,b_2)\in\R_+^2$, such that
\[{\sf CMP}_{(a_1,b_1)}\lesssim T_{\mu}(x)\lesssim{\sf CMP}_{(a_2,b_2)}.\]
Conversely, every CMP distribution is an element in $\cP^{1+}\cap\cP^{1-}$, and hence super-exponential. The Zeta distribution on $\N$ with parameter $a>1$ has probability mass function given by \cite{JKK05}:
\[{\sf Zeta}_a(x)=\frac{1}{\zeta(s)}x^{-a},\]
where $\zeta(a)=\sum_{i=1}^{\infty}i^{-a}$ is the Riemann zeta function of $a$.
For every probability distribution $\mu\in\cP^{3+}\cap\cP^{3-}$, there exists $a_1,\ a_2>1$ such that
\[{\sf Zeta}_{a_1}(x)\lesssim T_{\mu}(x)\lesssim{\sf Zeta}_{a_2}(x).\]
Conversely, every Zeta distribution is an element in $\cP^{3+}\cap\cP^{3-}$, and hence  sub-exponential.

\begin{figure}[h]
\begin{tikzpicture}[baseline=(current bounding box.north)]

\draw[fill=blue,opacity=.6] (0,0) -- (0,3) arc(90:270:3) --cycle;
\draw[fill=orange,opacity=.3] (0,0) -- (0,2) arc(90:270:2) --cycle;
\draw[fill=green,opacity=.3] (0,0) -- (0,1) arc(90:270:1) --cycle;

\draw[fill=yellow,opacity=.3]  (-3.6,0) -- (-3.6,3) arc(90:-90:3) --cycle;
\draw[fill=red,opacity=.3] (-3.6,0) -- (-3.6,2.5) arc(90:-90:2.5) --cycle;
\draw(-3.6,0) -- (-3.6,2) arc(90:-90:2) --cycle;
\draw[fill=brown,opacity=.6]  (-3.6,0) -- (-3.6,1.5) arc(90:-90:1.5) --cycle;


\node[] at (-.3,.7){$\cP^{1+}$};
\node[] at (-.6,1.5){$\cP^{2+}_{1}$};
\node[] at (-1.1,2.4){$\cP^{3+}$};

\node[] at (-3.2,.8){$\cP^{3-}$};
\node[] at (-3,1.6){$\cP^{2-}_{<1}$};
\node[] at (-2.9,2.1){$\cP^{2-}_{1}$};
\node[] at (-2.5,2.6){$\cP^{1-}$};

\node[] at (-.8,0){$I$};
\node[] at (-1.6,0){$II$};
\node[] at (-2.5,0){$III$};

\end{tikzpicture}
\caption{$I=\cP^{1+}\cap \cP^{1-}$: 
CMP-like distributions. $II=\cP^{2+}_{1}\cap \cP^{2-}_{1}$: Exponential-like distributions. $III=\cP^{3-}\cap \cP^{3+}$: Power-like distributions.}\label{fig_venn}
\end{figure}

We first provide the tail asymptotics for QSDs.
 
\thmb\label{th-19b}
Assume $\rm{(\mathbf{A1})}$-$\rm{(\mathbf{A3})}$ and  $\partial\neq\varnothing$. {Assume there exists a QSD $\nu$  of $Y_t$ on $\partial^{\sf c}$.\footnote{{Parameter conditions for the existence and ergodicity of QSDs are given  in \cite{XHW20a}}.}}  Then $\alpha\le0\le R$. Furthermore,  if
\begin{itemize}
\item[--] $R=0$, then $\alpha_-\ge\theta_{\nu}$. If, in addition, $R_-=R_+$, then  $\alpha\le-\theta_{\nu}$, and if $R_->R_+$, then $\beta\ge\theta_{\nu}$.
\item[--] $R_-=R_+>0$ and $\alpha=0$, then {$R\ge\sigman$}.
\end{itemize}
Moreover, if $R_->R_+$, and
\begin{enumerate}[label=(\roman*)]
\item $R=0$ and $\beta>\theta_{\nu}$, then $\nu\in\cP^{1-}_{(R_--R_+)\om_*^{-1}}$,
\item $R=0$, $\beta=\theta_{\nu}$ and $R_--R_+\le1$, then $\nu\in\cP^{2-}_{1}$,
\item $R=0$, $\beta=\theta_{\nu}$ and {$R_--R_+>1$, then $\nu\in\cP^{2-}_{R_--R_+-1}$,}
\item $R>0$, then $\nu\in\cP^{1-}_{(R_--R_+)\om_*^{-1}}$. If, in addition  $R>1$, then $\nu\in\cP^{1+}_{(R_--R_+)(\om_+\om_*)^{-1}}$.
\end{enumerate}
If $R_+=R_-$, and
\begin{enumerate}[label=(\roman*)]
\item[(v)] $R>0$ and $\alpha<0$, then $\nu\in\cP^{2-}_{1}$. If, in addition, $R>1$, then $\nu\in\cP^{2+}_{1}$,
\item[(vi)] $R>0$, $\alpha=0$, and
\begin{itemize}
\item[--] $R=\sigman<1$, then $\nu\in\cP^{2-}_{1-R}$,
\item[--] $R\ge\sigman=1$, then $\nu\in\cP^{3-}$,
\item[--] $\min\{1,R\}>\sigman$, then $\nu\in\cP^{2-}_{1<}$,
\end{itemize}
\item[(vii)] $R=0$, and
\begin{itemize}
\item[--] $\alpha+\theta_{\nu}=0$ and $\sigman<1$, then $\nu\in\cP^{2-}_{1-\sigman}$,
\item[--] $\alpha+\theta_{\nu}=0$ and  $\sigman=1$, then $\nu\in\cP^{3-}$,
\item[--] $\alpha+\theta_{\nu}<0$, then $\nu\in\cP^{2-}_{1}$.
\end{itemize}
\end{enumerate}
Furthermore, if
\begin{enumerate}[label=(\roman*)]
\item[(viii)] $R=1$, then $\nu\in\cP^{3+}_{\theta_{\nu}\alpha_+^{-1}}$,
\item[(ix)] $0<R<1$, then $\nu\in\cP^{2+}_{1-R}$,
\item[(x)] $R=0$ and $\alpha_->\theta_{\nu}$, then $\nu\in\cP^{2+}_{1}$,
\item[(xi)] $R=0$ and $\alpha_-=\theta_{\nu}$, then $\nu\in\cP^{1+}_{-\sigman\om_-^{-1}}$.
\end{enumerate}
\thme

\cob\label{co-3}
Assume $\rm{(\mathbf{A1})}$-$\rm{(\mathbf{A3})}$ and  $\partial\neq\varnothing$. No QSDs have a tail  decaying faster than a CMP distribution. {Moreover, any QSD, if such exists,} is super-exponential if $R_->\max\{1,R_+\}$ or (xi) holds, exponential-like if (a) $R_-=R_+=0$ and $\alpha+\theta_{\nu}<0$ or (b) $R_-=R_+>1$ and $\alpha<0$ holds,  and sub-exponential if (vi) or $R_-=R_+=\alpha+\theta_{\nu}=0$ holds, and in particular, it decays no faster than a power-like distribution if $R\ge \sigman=1$ and $\alpha=0$.
\coe

Analogously, we further characterise the tails of the stationary distributions.

\thmb\label{th-19}
Assume $\rm{(\mathbf{A1})}$-$\rm{(\mathbf{A2})}$ and $\partial=\varnothing$.  {Assume there exists a stationary distribution $\pi$  of $Y_t$ on $\cX$ with unbounded support.\footnote{{Parameter conditions for the existence and ergodicity of stationary distributions are given  in \cite{XHW20a}}.}} Then, $\alpha\le0$, and in particular, when $\alpha=0$,
\begin{itemize}
\item[--] if $\Delta=0$ then $\sigman<1$,
\item[--] if $\sigman=1$ then $\Delta>1$.
\end{itemize}
Moreover, if
\begin{enumerate}[label=(\roman*)]
\item $R_->R_+$, then $\pi\in\cP^{1+}_{(R_--R_+)(\om_+\om_*)^{-1}}\cap\cP^{1-}_{(R_--R_+)\om_*^{-1}}$,
\item $R_-=R_+$ and $\alpha<0$, then $\pi\in\cP^{2+}_{1}\cap\cP^{2-}_{1}$.
\item $\alpha=0$, $\Delta>0$, and $\sigman<1$,  then $\pi\in\cP^{2+}_{1-\sigman}\cap\cP^{2-}_{1-\sigman}$,
\item  $\alpha=0$ and $\sigman=1$, then $\pi\in\cP^{3+}_{\Delta-1}$. In particular, if in addition, (iv)' $\delta>1$, then $\pi\in\cP^{3+}_{\Delta-1}\cap\cP^{3-}_{\delta-1}$.
\item $\alpha=0$, $\Delta=0$, and $\rhon<1$, then $\pi\in\cP^{2-}_{1-\rhon}$.
\item $\alpha=0$, $\Delta=0$, and $\rhon\ge1$, then $\pi\in\cP^{3-}$.
\end{enumerate}
\thme
As a consequence, a trichotomy regarding the tails of the stationary distributions are derived.
\cob
Assume $\rm{(\mathbf{A1})}$-$\rm{(\mathbf{A2})}$ and $\partial=\varnothing$. Any stationary distribution of $Y_t$ on $\cX$ with unbounded  support, {if such exists}, is
\begin{itemize}
\item[--] super-exponential with a CMP-like tail if Theorem \ref{th-19}(i) holds,
\item[--] exponential-like if Theorem \ref{th-19}(ii) holds,
\item[--] sub-exponential if one of Theorem \ref{th-19}(iii), (iv)', and (vi) holds. In particular the tail is power-like if  (iv)' holds.
\end{itemize}
\coe
\cob\label{co-4}
Assume $\rm{(\mathbf{A1})}$, $\rm{(\mathbf{A2})'}$, $\rm{(\mathbf{A3})}$, $\partial=\varnothing$, $R\ge 3$, and $(R-1)\vartheta-\alpha_+\le0$. Any  stationary distribution of $Y_t$ on $\cX$ with unbounded  support, {if such exists}, is ergodic. {In particular, $Y_t$ is   non-explosive.}
\coe
\prb
By  $\rm{(\mathbf{A2})'}$, $R\in\N_0$ and $\sigman=1$. By \cite[Theorem~1]{XHW20a}, {under $\rm{(\mathbf{A3})}$,} $Y_t$ is explosive if either (1) $R\ge2$ and $\alpha>0$, or (2) $R\ge3$, $\alpha=0$ and $\gamma-\vartheta>0$.
By Theorem~\ref{th-19}, $\alpha\le0$. If a stationary distribution exists and if $Y_t$ is non-explosive, then  {under $\rm{(\mathbf{A3})}$,}  $Y_t$ is is positive recurrent and the stationary distribution is unique and ergodic \cite{N98}. When $\alpha=0$ and $R\ge3$, it follows from Theorem~\ref{th-19} that $\Delta>1$, i.e., $\gamma-\vartheta<(R-1)\vartheta-\alpha_+\le0$ as assumed. Hence, $Y_t$ is always non-explosive and thus the stationary distribution is ergodic.
\pre

We make the following remarks.


$\bullet$ The estimate of the tail does not depend on the choice of  ($R_{\om}^2,\ R_{\om}^3$) of the transition rate functions when $\alpha<0$, whereas it may depend when $\alpha=0$. In this case, the larger $\sigman$ and $\rhon$ are, the sharper the results are.

$\bullet$ Generically, no limit distributions (in the cases covered) can decay faster than a CMP distribution nor slower than a Zeta distribution.

$\bullet$ The unique gap case in Theorem~\ref{th-19} is $\alpha=0$, $\sigman<1$ and $\Delta=0$.

$\bullet$ Assume ($\rm\mathbf{A2}$)'. By Corollary~\ref{co-4}, if the chain $Y_t$ is explosive and a stationary distribution exists, then $\alpha=0$ and $R\ge3$.

$\bullet$ Although not stated explicitly in Theorem~\ref{th-19}, the tail asymptotics of a stationary distribution of a BDP (cases (i)-(iii) and (vi)') is sharp up to the leading order, in comparison with Proposition~\ref{pro-15}. Similarly, when $R_->R_+$, then the tail asymptotics is sharp up to the leading order for {\em upwardly skip-free processes} (case (i)) \cite{A91}. In comparison with the sharp results provided in Proposition~\ref{pro-16}, the results  obtained in Theorem~\ref{th-19}   capture the leading tail asymptotics  of a stationary distribution, e.g., in case (iii).

$\bullet$ The assumption that $R>1$ in Corollary~\ref{co-3} is crucial. Indeed, as Examples~\ref{ex-1} and \ref{ex-2} illustrate, when $R=1$ and $\alpha<0$, the QSD may still exist and has either geometric or Zeta-like tail. This means $\alpha<0$ is not sufficient for any QSD to have an exponential tail. It remains  to see if a QSD with CMP tail may exist when $R=1$. Moreover, we emphasise that $R>1$ and $\alpha<0$ ensure the existence of a unique ergodic QSD assuming ($\rm\mathbf{A2}$)' \cite{XHW20b}. Hence, such ergodic QSDs are not heavy-tailed.

$\bullet$  Let $Y_t$ be a CTMC and $\widetilde{Y}_t$ a BDP on the same state space. If the critical parameters are the same ($\alpha, \beta,  R,\ldots$), then the two processes share the same qualitative characterisation in terms of existence of limit distributions, ergodicity, explosivity, \ldots \cite{XHW20b}. One might conjecture that the decay of limit distributions (provided such exist) is also the same and takes the sharp form induced by the BDP. In other words, one might conjecture that the tail asymptotics of a general CTMC coincide with that of a BDP \emph{representative}. 
 
In the following, we illustrate and elaborate on the results by example.
 The examples have real analytic APLH transition rate functions, and thus $\sigman=1$ and $\rhon=2$.

\smallskip

Assumption $\rm{(\mathbf{A1})}$  is crucial for Theorem~\ref{th-19}.
\eb
Consider a model of stochastic gene regulatory expression \cite{SS08}, given by $\cT=\{-1\}\cup\N$ and
$$q(x,x-1)=\mathbbm{1}_{\N}(x),\quad q(x,x+j)=ab_{j-1},\quad  j\in\N,\ x\in\N_0,$$
where $a>0$, and $b_j\ge0$ for all $j\in\N_0$. Here the backward jump $-1$ represents the degradation of mRNA with unity degradation rate, and the forward jumps $j\in\N$ account for bursty production of mRNA  with transcription rate $a$ and burst size distribution $(b_j)_{j\in\N_0}$.  When $b_j=(1-\delta)\delta^j$ for  $0<\delta<1$, the stationary distribution is the {\em negative binomial distribution} \cite{SS08}:
\[ \pi(x)=\frac{\Gamma(x+a)}{\Gamma(x+1)\Gamma(a)}\delta^x(1-\delta)^a,\quad x\in\N_0.\]
When $a=1$, $\pi$ is also geometric. While Theorem~\ref{th-19}, if it did apply, would seem to suggest $T_{\pi}$ decays like a CMP distribution, since $1=R_->R_+=0$. The technical reason behinds this, is that  the proof  applies Corollary~\ref{co-5} which requires $\rm{(\mathbf{A1})}$.
It does not seem possible to directly extend the result in Theorem~\ref{th-19} to  CTMCs with unbounded jumps.
\ee

\eb
Consider a BDP with birth and death rates:
$$q(x,x-1)=\sum_{j=1}^bS(b,j)x^{\underline{j}},\quad q(x,x+1)=a,\quad x\in\N_0,$$
where  $a>0$, $b\in\N$, and $S(i,j)$ is the Stirling numbers of the second kind \cite{AS72}. Here $R_+=0<R_-=b$, and $\alpha=-S(b,b)=-1<0$. 
This BDP has an ergodic stationary distribution on $\N_0$ \cite{XHW20b}, and the unique stationary distribution is $\pi={\sf CMP}_{a,b}$.
\ee

By  Theorem~\ref{th-19b}(v), the tail of a QSD decays no faster than exponential distributions when $\alpha<0$ and $0\le R_-=R_+\le1$, which is also confirmed by the examples below.
\eb\label{ex-1}
Consider the linear BDP on $\N_0$ with birth and death rates:
$$q(x,x+1)=b x,\quad x\in\N_0,\quad \text{and}\quad q(1,0)=d,\quad q(x,x-1)=\lt(d\cdot2^{-1}+b\rt)(x+1),\quad x\in\N\!\setminus\!\!\{1\},$$
where $b$ and $d$ are positive constants \cite{O07}. For this process, a QSD $\nu$ is
\[\nu(x)=\frac{1}{x(x+1)},\quad x\in\N.\]
Hence $T_{\nu}$ decays as fast as the Zeta distribution with parameter $2$. Here $\alpha=-d/2<0$, and $R=R_-=R_+=1$.
\ee

\eb  \label{ex-2}
 Consider the linear BDP on $\N_0$ with $b_j=b j$ and $d_j=d_1 j$ with $0<b<d_1$ and $j\in\N$ \cite{V91}. For this process, a QSD $\nu$ is
  \[\nu(x)=\lt(\frac{b}{d_1}\rt)^{\!\!x-1}\lt(1-\frac{b}{d_1}\rt),\quad x\in\N,\]
 a geometric distribution. Here $\alpha=b-d_1<0$, and $R=R_-=R_+=1$.
\ee
By Theorem~\ref{th-19b}(iv) and (viii), the tail of the QSD decays no faster than CMP distributions and no more slowly than Zeta distribution, if $R_+<R_-=1$, which is also confirmed by the example below.
\eb
Consider a BDP on $\N_0$:
$$q(x,x+1)=\frac{x}{x+2},\quad x\in\N_0;\quad q(x,x-1)=x-1+2\frac{1}{x},\quad x\in\N.$$
Here $R=R_-=1>R_+=0$, $\alpha=-1$, $\sigman=1$, $\rhon=2$. Using the same Lyapunov function constructed in the proof of \cite[Theorem~4.4]{XHW20b}, it can be shown that there exists a uniquely ergodic QSD. By Theorem~\ref{th-19b}(iv), the tail of the QSD decays no more slowly than a CMP distribution. Indeed, the QSD is given by
\[\nu(x)=\frac{1}{(x-1)!(x+1)},\quad T_{\nu}(x)=\frac{1}{x!},\quad x\in\N.\]
The tail of the QSD decays like a Poisson distribution.
\ee

\eb
Consider a BDP on $\N_0$:
$$q(x,x+1)=x^2,\quad q(x,x-1)=x^2+x,\quad x\in\N_0.$$  Here $R=R_-=R_+=2>1$, and $\alpha=0$. Corollary~\ref{co-3} states that any QSD (if it exists) is heavy-tailed and its tail decays no faster than a Zeta-like distribution. Indeed, a QSD of the process is given by
\[\nu(x)=\frac{1}{x(x+1)},\quad T_{\nu}(x)=\frac{1}{x},\quad x\in\N.\]

\ee

\eb
Consider a quadratic BDP on $\N_0$:
$$q(x,x+1)=x(x+3)/2,\quad q(x,x-1)=x(x+1),\quad x\in\N_0.$$ Then $R_-=R_+=R=2>1$, $\alpha=-1/2$. Hence there exists a uniquely ergodic QSD \cite{XHW20b}. By Corollary~\ref{co-3}, this QSD decays exponentially. Indeed, the QSD is given by
\[\nu(x)=2^{-x},\quad T_{\nu}(x)=2^{-x+1},\quad x\in\N.\]
\ee

\section{Applications}\label{sect-appl}

In this section, we apply the results on  asymptotic tails to diverse models in biology. We emphasise that for all models/applications, the transition rate functions are real analytic APLH on a subset of $\N_0$, and thus $\sigman=1$ and $\rhon=2$, which we will not further  mention explicitly.

\subsection{Biochemical reaction networks}

In this section, we apply the results of Section~\ref{sect-Asymptotics} to some examples of  {\em stochastic reaction networks} (SRNs) with mass-action kinetics. These  are use to describe interactions of constituent molecular species with many applications in systems biology, biochemistry, genetics and beyond  \cite{G83,PCMV15}. An SRN with mass-action kinetics is a CTMC on $\N_0^d$ ($d\ge 1$) encoded  by a labeled directed graph
 \cite{AK11}.  We concentrate on SRNs on $\N_0$ with one species (\text{S}). In this case the graph is composed of reactions (edges) of the form  $n\text{S}\ce{->[\kappa]} m\text{S}$, $n,m\in\N_0$ ($n$ molecules of species \text{S} is converted into $m$ molecules of the same species),  encoding a jump from $x$ to $x+m-n$ with propensity  $q(x,x+m-n)=\kappa x(x-1)\ldots(x-n+1)$, $\kappa>0$. Note that multiple reactions might result in the same jump vector.

In general little is  known about the stationary distributions of  a reaction network, let alone the QSDs, provided either such exist \cite{ACK10,GBK14,HW20,HM19}. Special cases  include \emph{complex balanced} networks (in arbitrary dimension) which have Poisson product-form distributions \cite{ACK10,CW16}, reaction networks that are also birth-death processes, and reaction networks with  irreducible components, each with a finite number of states.

\eb\label{Ex-WC}
To show how general the results are we consider two SRNs, none of which are birth-death processes.
{\rm(i)} Consider a reaction network with a  strongly connected graph \cite{XHW20b}:
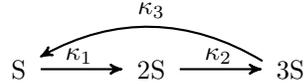
\[
\begin{tikzpicture}[node distance=3.5em, auto, scale=1]
 \tikzset{
    >=stealth',
    pil/.style={
           ->,
           thick,
           shorten <=2pt,
           shorten >=2pt,}
}
 \node[] (a) {};
  \node[right=1.3em of a] (n1) {};
  \node[above=-.5em of n1] (m1) {$\ka_1$};
\node[right=of a] (b) {2S};
\node[right=1.3em of b] (n2) {};
  \node[above=-.5em of n2] (m2) {$\ka_2$};
  \node[above=.9em of b] (m3) {$\ka_3$};
  \node[left=of b] (aa) {S} edge[pil, black, bend left=0] (b);
\node[right=of b] (c) {3S};
 \node[left=of c] (bb) {} edge[pil, black, bend left=0] (c);
\node[right=of b] (cc) {} edge[pil, black, bend right=30] (aa);
\end{tikzpicture}\]

For this reaction network, $\cT=\{1,-2\}$, and
\[q(x,x+1)=\kappa_1 x+\kappa_2 x(x-1),\quad q(x,x-2)=\kappa_3 x(x-1)(x-2).\]
Hence, $\alpha=-\kappa_3<0$. It is known that there exists s a unique exponentially ergodic stationary distribution $\pi$ on $\N$ \cite{XHW20b}. (The state $0$ is neutral \cite{XHW20a}.) By Theorem~\ref{th-19}, $\pi\in\cP^{1+}_1\cap\cP^{1-}_1$. Hence $\pi$ is light-tailed and $T_{\pi}$ decays as fast as a Poisson distribution.
since $\om_+=\om_*=1$, $R_+=2$, $R_-=3$. However, the stationary distribution is generally {\em not} Poisson. If $\kappa_2^2=\kappa_1\kappa_3$, then the reaction network  is \emph{complex-balanced}, hence the stationary distribution is Poisson \cite{ACK10}. If the parameter identity is not fulfilled then the distribution cannot be Poisson in this case \cite{CW16}.

{\rm(ii)} Consider a similar reaction network including direct degradation of S \cite{XHW20b}:
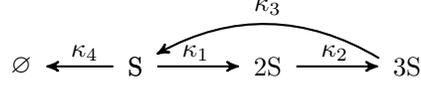
\[
\begin{tikzpicture}[node distance=3.5em, auto, scale=1]
 \tikzset{
    >=stealth',
    pil/.style={
           ->,
           thick,
           shorten <=2pt,
           shorten >=2pt,}
}
\node[] (a) {};
\node[right=1.3em of a] (n1) {};
\node[left=of a] (d) {$\varnothing$};
\node[right=1.3em of d] (n4) {};
\node[above=-.5em of n4] (m4) {$\ka_{4}$};
\node[right=1.3em of b] (n2) {};
\node[above=-.5em of n1] (m1) {$\ka_1$};
\node[right=of a] (b) {2S};
\node[right=1.3em of b] (n2) {};
\node[above=-.5em of n2] (m2) {$\ka_2$};
\node[above=.9em of b] (m3) {$\ka_3$}; \node[left=of b] (aa) {S} edge[pil, black, bend left=0] (d);
  \node[left=of b] (aa) {S} edge[pil, black, bend left=0] (b);
\node[right=of b] (c) {3S};
 \node[left=of c] (bb) {} edge[pil, black, bend left=0] (c);
\node[right=of b] (cc) {} edge[pil, black, bend right=30] (aa);
\end{tikzpicture}\]
The threshold parameters are the same as in {\rm(i)}, and it follows from \cite{XHW20b} that the reaction network has a uniformly exponentially ergodic QSD $\nu$. By Theorem \ref{th-19b}, $T_{\nu}$ decays like a CMP distribution.
\ee

\begin{example}
The following bursty Schl\"{o}gl model was proposed in \cite{FMD17}:
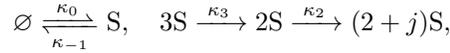
\begin{equation*}
\varnothing\ce{<=>[\ka_0][\ka_{-1}]}\tS,\quad 3\tS\ce{->[\ka_3]}2\tS\ce{->[\ka_2]}(2+j)\tS,
\end{equation*}
where $j\in\N$. When $j=1$, it reduces to the classical Schl\"{o}gl model.

The associated process has  a unique ergodic stationary distribution  $\pi$ on $\N_0$ \cite{XHW20b}. Bifurcation with respect to patterns of the ergodic stationary distribution is discussed in \cite{FMD17}, based on a diffusion approximation in terms of the Fokker-Planck equation.  Using the results established in this paper the tail distribution can be characterised rigorously. In fact $\pi\in\cP^{1+}_{j^{-1}}\cap\cP^{1-}_1$.  Hence, $\pi$ is light-tailed and $T_{\pi}$ decays like a CMP distribution.

\prb
We have $\cT=\{-1,1\}\cup\{j\}$. The ergodicity follows from \cite[subsection~4.3]{XHW20b} as a special case. It is straightforward  to verify that  ${\rm(\mathbf{A1})}$-${\rm(\mathbf{A3})}$ are all satisfied. Moreover, $R_+=2$, $R_-=3$. Then the conclusion follows from Theorem~\ref{th-19}.
\pre
\end{example}

\begin{example}
{Consider the following one-species SRN:
\[\tS\ce{->[$\kappa_1$]}4\tS,\quad 3\tS\ce{->[$\kappa_2$]}0.\]
In this case, $\omega_*=3$  and the ambient space $\N_0$ is divided into disjoint irreducible sets, $3\mathbb{N}$, $3\mathbb{N}_0+1$ and  $3\mathbb{N}_0+2$, as well as an absorbing state $0$. 
 By simple calculation, we have $R_-=3>R_+=1$, and $\alpha = -3\kappa_2<0$. Hence, this SRN is   positive recurrent on the sets  $3\mathbb{N}_0+1$  and  $3\mathbb{N}_0+2$, while it admits a unique QSD on   $3\mathbb{N}$ \cite[Theorem~7]{XHW20a}. According to Theorem \ref{th-19}, the tails of the stationary distributions decay like CMP on $3\mathbb{N}_0+1$  and  $3\mathbb{N}_0+2$, and according to Theorem \ref{th-19b}, the tail of the QSD also decays as CMP on  $3\mathbb{N}$. Since the transition rate functions take a common form (polynomial) on all three  irreducible sets, then the parameters are the same on all three irreducible sets (open and closed), and consequently, the tail asymptotics are the same on the two closed irreducible sets. The same would hold true for the open irreducible sets, if there were more than one open irreducible set.}
\end{example}

\eb
Consider the following one-species  S-system modelling a gene regulatory network \cite{CCE15}:
\[\varnothing\ce{->[(\ka_1,\xi_1)]}\tS\ce{<-[(\ka_2,\xi_2)]}3\tS\]
with the following \emph{generalised mass action kinetics} (GMAK):
$$q(x,x+1)=\ka_1\frac{\Gamma(x+\xi_1)}{\Gamma(x)},\quad q(x,x-2)=\ka_2\frac{\Gamma(x+\xi_2)}{\Gamma(x)},\quad x\in\N_0,$$
where $\ka_1, \ka_2>0$ are the reaction rate constants, and $\xi_2>\xi_1>0$ are the indices of GMAK.

By Stirling's formula
\[\log\Gamma(x)=(x-1/2)\log x-x+\log\sqrt{2\pi}+\rO(x^{-1}),\]
hence $q(x,x+1)$ is APLH with $(\xi_1,\xi_1-1,\xi_1-2)$ and $q(x,x-2)$ is APLH with $(\xi_2,\xi_2-1,\xi_2-2)$.  Then $R_-=\xi_2>R_+=\xi_1$, $\om_*=1$, $\om_-=-2$ and $\om_+=1$. Using the same Lyapunov function constructed in the proof of \cite[Theorem~4.4]{XHW20b}, it can be shown that there exists a uniquely ergodic stationary distribution $\pi$ on $\N_0$ with support $\N$. By Theorem~\ref{th-19}, $\pi\in\cP^{1+}_{\xi_3-\xi_1}\cap\cP^{1-}_{\xi_3-\xi_1}$.
\ee

\subsection{An extended class of branching processes}

Consider an extended class of branching processes on $\N_0$ \cite{C97} with transition rate matrix $Q=(q(x,y))_{x,y\in\N_0}$:
\[q(x,y)=\left\{\begin{array}{cl}
r(x)\mu(y-x+1), &\quad  \text{if}\quad  y\ge x-1\ge0\quad \text{and}\quad y\neq x,\\
 -r(x)(1-\mu(1)), &\quad  \text{if}\quad y=x\ge1,\\
q(0,y), & \quad \text{if}\quad y>x=0,\\
 -q(0), &\quad \text{if}\quad y=x=0,\\
 0, &\quad  \text{otherwise},
\end{array}\right.\]
where $\mu$ is a probability measure on $\N_0$, $q(0)=\sum_{y\in\N}q(0,y)$, and $r(x)$ is a positive finite function on $\N_0$.
Assume

\medskip
\noindent($\rm\mathbf{H1}$) $\mu(0)>0$, $\mu(0)+\mu(1)<1$.

\medskip
\noindent($\rm\mathbf{H2}$) $\sum_{y\in\N}q(0,y)y<\infty$, $M
=\sum_{k\in\N_0}k\mu(k)<\infty$.
\medskip

\noindent($\rm\mathbf{H3}$) $r(x)$ is a polynomial of degree $R\ge1$ for large $x$.
\medskip

The tail asymptotics  of infinite stationary measures in the null recurrent case is investigated in \cite{LZ11} under ($\rm\mathbf{H1}$)-($\rm\mathbf{H2}$) for general $r$. Here we assume $r$ is polynomial ($\rm\mathbf{H3}$).
The  following is a consequence of the results of Section~5.

\thmb
Assume ${\rm(\mathbf{H1})}$-${\rm(\mathbf{H3})}$, $Y_0\neq0$, and that $\mu$ has finite support. \enb
\item[(i)] Assume $q_0>0$. Then there exists an ergodic stationary distribution $\pi$ on $\N_0$ if (i-1) $M<1$ or (i-2) $M=1$ and $R>1$. Moreover, $T_{\pi}$ decays like a geometric distribution if (i-1) holds while like a Zeta distribution if (i-2) holds.
\item[(ii)] Assume $q_0=0$. Then there exists an ergodic QSD $\nu$ on $\N$ if (i-1) $M<1$ and $R>1$ or (i-2) $M=1$ and $R>2$. Moreover, $T_{\nu}$ decays like a geometric distribution if (ii-1) holds while no faster than a Zeta distribution if (ii-2) holds.
\ene
\thme

\prb
For all $k\in\cT$, let
$$q(x,x+k)=\cab r(x)\mu(k+1),\quad \text{if}\ x\in\N,\\ q(0,k),\qquad\qquad\quad\, \text{if}\ x=0.\cae$$
 By ${\rm(\mathbf{H1})}$, $\mu(k)>0$ for some $k\in\N$. Hence regardless of $q(0)$, by positivity of $r$, ${\rm(\mathbf{A1})}$-${\rm(\mathbf{A3})}$ are  satisfied with $\cT_-=\{-1\}$ and $\cT_+=\{j\in\N\colon j+1\in\supp\mu\,\,\, \text{or}\,\,\, q(0,j)>0\}$.
  Let $r(x)=ax^R+bx^{R-1}+\rO(x^{R-2})$ with $a>0$. It is straightforward to verify that
$R_+=R_-=R,$ $\alpha=a(M-1)$. The ergodicity follows from \cite{XHW20b}, and the tail asymptotics follow from Theorems~\ref{th-19} and \ref{th-19b}.
\pre

\subsection{Stochastic population processes under bursty reproduction}
Two stochastic population models with bursty reproduction are investigated in \cite{BA16}.

The first model is a Verhulst logistic population process with bursty reproduction. The process $Y_t$ is a CTMC on $\N_0$ with transition rate matrix $Q=(q(x,y))_{x,y\in\N_0}$ satisfying:
\[q(x,y)=\cab c\mu(j)x,\qquad \text{if}\ y=x+j,\ j\in\N,\\ \frac{c}{K}x^2+x,\quad\, \text{if}\ y=x-1\in\N_0,\\ 0,\qquad\qquad\ \, \text{otherwise},\cae\]
where $c>0$ is the  reproduction rate, $K\in\N$ is the typical population size in the long-lived
metastable state prior to extinction \cite{BA16}, and $\mu$ is the burst size  distribution.

Approximations of the mean time to extinction and QSD are discussed in \cite{BA16} against various different burst size distributions  of finite mean (e.g., Dirac measure, Poisson distribution, geometric distribution, negative-binomial distribution). The existence of an ergodic QSD for this population model is established in \cite{XHW20b}. Nevertheless, the tails of QSD is not addressed therein.

\thmb
Assume $\mu$ has a finite support. Let $\nu$ be the unique ergodic QSD on $\N$ trapped to zero for the Verhulst logistic model $Y_t$. Then $T_{\nu}$ decays like a CMP distribution.\thme
\prb
We have $\cT=\{-1\}\cup\supp\mu$, $q(x,x-1)=\frac{c}{K}x^2+x$, $q(x,x+k)=c\mu(k)x$, for $k\in\supp\mu$ and $x\in\N$. Since $\mu$ has a finite support, ${\rm(\mathbf{A1})}$-${\rm(\mathbf{A3})}$ are satisfied. Moreover, since $\supp\mu\neq\varnothing$, we have $R_-=2$ and $R_+=1$. Again, the ergodicity result follows from \cite{XHW20b}. The tail asymptotics follow directly from  Theorem~\ref{th-19}.
\pre

In subsequent sections, we provide proofs of the main results in Section~\ref{sect-Asymptotics}. Since the proof of Theorem~\ref{th-19b} is based on that of Theorem~\ref{th-19}. We first prove Theorem~\ref{th-19} in the next section.

\section{Proof of Theorem~\ref{th-19}}

Let 
\[\begin{split}
  \alpha_j=&\cab\lim_{x\to\infty}\frac{\sum_{\om\in A_j}q(x,x+\omega)}{x^{R_-}},\quad \text{if}\  j=\omega_-+1,\ldots,0,\\
  \lim_{x\to\infty}\frac{\sum_{\om\in A_j}q(x,x+\omega)}{x^{R_+}},\quad \text{if}\  j=1,\ldots,\omega_+,\\ 0,\qquad\qquad\qquad\qquad\qquad\quad \text{otherwise},
  \cae\\
  \gamma_j=&\cab\lim_{x\to\infty}\frac{\sum_{\om\in A_j}q(x,x+\omega)-\alpha_j x^{R_-}}{x^{R_--\sigman}},\quad \text{if}\ j=\omega_-+1,\ldots,0,\\ \lim_{x\to\infty}\frac{\sum_{\om\in A_j}q(x,x+\omega)-\alpha_j x^{R_+}}{x^{R_+-\sigman}},\quad \text{if}\ j=1,\ldots,\omega_+,\\ 0,\hspace{4.8cm} \ \text{otherwise}.\cae
\end{split}\]

Note that $\beta=\alpha_0$.
From Lemma~\ref{Sle-1}, $\alpha_-=\om_*\sum_{j=\omega_-+1}^0\alpha_{j}$, $\alpha_+=\om_*\sum_{j=1}^{\omega_+}\alpha_j$. By ${\rm(\mathbf{A3})}$, \eqb\label{Eq-3}\sum_{\om\in A_j}q(x,x+\omega)=\cab x^{R_-}(\alpha_j+\gamma_jx^{-\sigman}+\rO(x^{-\rhon})),\quad \text{if}\ j=\om_-+1,\ldots,0,\\ x^{R_+}(\alpha_j+\gamma_jx^{-\sigman}+\rO(x^{-\rhon})),\quad\, \text{if}\ j=1,\ldots,\om_+.\cae\eqe
Since $$q(x,x+j\omega_*)={\rm sgn}(j)\Bigl(\sum_{\om\in A_j}q(x,x+\omega)-\sum_{\om\in A_{j+1}}q(x,x+\omega)\Bigr),\quad j=\om_-,\ldots,-1,1,\ldots,\om_+,$$ we have
$$q(x,x+j\omega_*)=\cab x^{R_-}((\alpha_{j+1}-\alpha_j)+(\gamma_{j+1}-\gamma_j)x^{-\sigman}+\rO(x^{-\rhon})),\quad \text{if}\ j=\om_-,\ldots,-1,\\ x^{R_+}((\alpha_j-\alpha_{j+1})+(\gamma_j-\gamma_{j+1})x^{-\sigman}+\rO(x^{-\rhon})),\quad\, \text{if}\ j=1,\ldots,\om_+.\cae$$

Since ($\rm\mathbf{A1}$)-($\rm\mathbf{A2}$) imply ${\cap}_{\om\in\cT}\ \{x\in\cX\colon q(x,x+\omega)=0\}$ is finite, {then 
 it easily follows} that both $\cT_-\neq\varnothing$ and $\cT_+\neq\varnothing$ since $\supp\pi$ is unbounded. Hence $\alpha_-\ge\alpha_0>0$, $\alpha_+\ge\alpha_1>0$, and $-\infty<\om_-<\om_+<\infty$.

For the ease of exposition and w.o.l.g., we assume throughout the proof that $\om_*=1$ (recall the argument in the proof of Theorem~\ref{th-18}). Hence $\N_0+b\subseteq\cX\subseteq\N_0$ for some $b\in\N_0$ by Proposition~\ref{pro-0}.

Most inequalities below are based on the identities in Theorem~\ref{th-18} and Corollary~\ref{co-5}. Therefore, we use LHS (RHS) with a label in the subscript as shorthand for the left (right) hand side of an equation with the given label.

The claims that $\Delta=0$ implies $\sigman<1$ and that $\sigman=1$ implies $\Delta>1$ are proved in (iii)-(vi) Step I below.

We first show $\alpha\le0$. Suppose by means of contradiction that $\alpha>0$. Then either (1) $R_+>R_-$ or (2) $R_+=R_-$ and $\alpha_+>\alpha_-$ holds.
Define the auxiliary function
\[f_j(x)=\sum_{\om\in A_j}q(x-j,x-j+\omega),\quad j=\om_-+1,\ldots,\om_+.\]
From \eqref{Eq-3} it follows that
\[
f_j(x)=\cab x^{R_-}(\alpha_j+\gamma_jx^{-\sigman}-\alpha_jj{R_-}x^{-1}+\rO(x^{-\min\{\rhon,\sigman+1\}})),\quad \text{if}\ j=\om_-+1,\ldots,0,\\ x^{R_+}(\alpha_j+\gamma_jx^{-\sigman}-\alpha_jj{R_+} x^{-1}+\rO(x^{-\min\{\rhon,\sigman+1\}})),\quad \text{if}\ j=1,\ldots,\om_+.\cae\]
Let
$$\beta_j(x)=\cab x^{-R_-}f_j(x)-\alpha_j,\quad \text{if}\ j=\om_-+1,\ldots,0,\\ x^{-R_+}f_j(x)-\alpha_j,\quad \text{if}\ j=1,\ldots,\om_+.\cae$$
 Then, there exist $N_3,\ N_4\in\N$ with $N_3>N_1, N_4$ such that for all $x\ge N_3$,
\eqb\label{Eq-38}|\beta_j(x)|\le N_4x^{-\sigman},\quad  j=\om_-+1,\ldots,\om_+,\eqe
From {Theorem~\ref{th-18}{\rm(3)}}, we have
 \eqb\label{Eq-46-a}x^{R_--R_+}\sum_{j=\omega_-+1}^0\lt(\alpha_{j}+\beta_{j}(x)\rt)\pi\lt(x-j\rt)=
\sum_{j=1}^{\omega_+}\lt(\alpha_j+\beta_j(x)\rt)\pi\lt(x-j\rt),\eqe
Since $R_-\le R_+$ and $T_{\pi}(x)\le1$ for all $x\in\N_0$, summing up in \eqref{Eq-46-a} from $x$ to infinity yields \eqb\label{Eq-4}
\sum_{y=x}^{\infty}y^{R_--R_+}\sum_{j=\omega_-+1}^0\lt(\alpha_{j}+\beta_{j}(y)\rt)\pi\lt(y-j\rt)
=\sum_{y=x}^{\infty}
\sum_{j=1}^{\omega_+}\lt(\alpha_j+\beta_j(y)\rt)\pi\lt(y-j\rt).\eqe
In the light of the monotonicity of $T_{\pi}(x)$ and $x^{R_--R_+}$, it follows from \eqref{Eq-38} that there exists $C=C(N_4)>0$ and $N_5\in\N$ with $N_5\ge N_3$ such that for all $x\ge N_5$,
\[\begin{split}
\text{LHS}_{\eqref{Eq-4}}& \le\sum_{y=x}^{\infty}y^{R_--R_+}\sum_{j=\omega_-+1}^0\lt(\alpha_j+N_4y^{-\sigman}\rt)\pi\lt(y-j\rt)\\
&  \le x^{R_--R_+}\sum_{j=\omega_-+1}^0\lt(\alpha_j+N_4x^{-\sigman}\rt)\sum_{y=x}^{\infty}\pi\lt(y-j\rt)\\
& =x^{R_--R_+}\sum_{j=\omega_-+1}^0\lt(\alpha_j+N_4x^{-\sigman}\rt)T_{\pi}(x-j)\\
& \le x^{R_--R_+}T_{\pi}(x)\sum_{j=\omega_-+1}^0\lt(\alpha_j+N_4x^{-\sigman}\rt)\\
& \le x^{R_--R_+}T_{\pi}(x)\lt(\alpha_-+Cx^{-\sigman}\rt).\end{split}\]
Similarly, with a possibly larger $C$ and $N_5$, for all $x\ge N_5$, {one can show
\[\begin{split}
\text{RHS}_{\eqref{Eq-4}}
&\ge\lt(\alpha_+-Cx^{-\sigman}\rt)T_{\pi}(x-1),
\end{split}\]}
which further implies that for all $x$ large enough,
\[1\ge\frac{T_{\pi}(x)}{T_{\pi}(x-1)}\ge x^{R_+-R_-}\frac{\alpha_+-Cx^{-\sigman}}{\alpha_-+Cx^{-\sigman}}>1,\]
since either (1) $R_+>R_-$ or (2) $R_+=R_-$ and $\alpha_+>\alpha_-$ holds.  This  contradiction shows that $\alpha\le0$.

\medskip

Next, we provide asymptotics of $T_{\pi}(x)$ case by case.
\smallskip
{
Beforehand, let us illustrate the   idea behind the proof by means of an l example: the BDP case ($\omega_+=-\omega_-=1$) for $R_->R_+$. Again, assume w.l.o.g. that $\omega_*=1$. From Corollary~\ref{co-5} it follows that 
\[\frac{T_{\pi}(x)(q(x,x-1)+ q(x-1,x))}{T_{\pi}(x-1)q(x-1,x)} - \frac{T_{\pi}(x+1) q(x,x-1)} {T_{\pi}(x-1)q(x-1,x)} = 1,\]
which implies by non-negativity of $T_{\pi}(x+1) q(x,x-1)$ that 
\[\frac{T_{\pi}(x)}{T_{\pi}(x-1)} \ge \frac{q(x-1,x)}{q(x-1,x) + q(x,x-1)}.\]
Note that 
$$q(x-1,x) = (x-1)^{R_+} (\alpha_1 +\gamma_1 x^{-\sigma_1} + O(x^{-\sigma_2})),$$
$$q(x,x-1) = x^{R_-}(\alpha_0+\gamma_0 x^{-\sigma_1}+O(x^{-\sigma_2})).$$
This further shows that 
\[T_{\pi}(x) \gtrsim \Gamma(x)^{R_+-R_-}\lt(\frac{\alpha_1}{\alpha_0}\rt)^{x+Cx^{1-(R_--R_+)}+\rO(\log x)}\] 
 for some constant $C>0$.
}

{
To obtain the upper estimate, rewrite \eqref{Eq-46-a} as \[(\alpha_0+\beta_0(x))\pi(x) = (\alpha_1+\beta_1(x))x^{R_+-R_-}\pi(x-1)\] 
Summing up the above equation from $x$ to infinity yields that 
\[\sum_{y\ge x}(\alpha_0+\beta_0(y))\pi(y) = \sum_{y\ge x} (\alpha_1+\beta_1(x))x^{R_+-R_-}\pi(y-1).\]
From \eqref{Eq-38}, it follows that there exists $C_1,\ N>0$ such that for all large $x\ge N$,
\begin{align*}
\sum_{y\ge x}(\alpha_0+\beta_0(y))\pi(y) \ge \lt(\alpha_0-C_1x^{-\sigman}\rt)T_{\pi}(x),
\end{align*}
while 
\begin{align*}
\sum_{y\ge x} (\alpha_1+\beta_1(x))x^{R_+-R_-}\pi(y-1) \le& \sum_{y=x}^{\infty}y^{R_+-R_-}\lt(\alpha_1+C_1y^{-\sigman}\rt)\pi\lt(y-1\rt)\\
& \le x^{R_+-R_-}T_{\pi}(x-1)\lt(\alpha_1+C_1x^{-\sigman}\rt)
\end{align*}
Hence,
\[\frac{T_{\pi}(x)}{T_{\pi}(x-1)}\le x^{R_+-R_-}\frac{\alpha_1+C_1x^{-\sigman}}{\alpha_{0}-C_1x^{-\sigman}}=x^{R_+-R_-}\lt(\frac{\alpha_1}{\alpha_{0}}+
\rO(x^{-\sigman})\rt),\]
which further implies that 
\[T_{\pi}(x)(x)\lesssim\Gamma(x)^{R_+-R_-}\lt(\frac{\alpha_1}{\alpha_0}\rt)^{x
+C_2x^{1-\sigman}+\rO(\log x)}\] for some $C_2\ge C_1$.
}
\medskip

Now we provide the detailed asymptotic estimates of the tail distribution case by case.

\noindent{\rm(i)} $R_->R_+$. Recall  Stirling's formula for the Gamma function \cite{AS72}:
\[\log\Gamma(x)=x\log x-x+\rO(\log x),\] where $\log$ is the natural logarithm.
Based on this Stirling's formula, it suffices to prove that there exists $\widetilde{C}>0$ such that
 \eqb\label{Eq-27-a}
T_{\pi}(x)\gtrsim\Gamma(x)^{R_+-R_-}\lt(\frac{\alpha_1}{\alpha_0}\rt)^{x+\widetilde{C}x^{1-(R_--R_+)}+\rO(\log x)},\eqe
\eqb\label{Eq-27-b}T_{\pi}(x)\lesssim\Gamma(x\om_+^{-1})^{R_+-R_-}\lt(\frac{\alpha_+}{\alpha_0}\omega_+^{R_+-R_-}\rt)^{\om^{-1}_+x
+\widetilde{C}x^{1-\sigman}+\rO(\log x)}.
\eqe
Next, we prove \eqref{Eq-27-a} and \eqref{Eq-27-b} one by one.

We first show \eqref{Eq-27-a}. Recall that ($\rm\mathbf{A2}$) ensures that there exists $N\in\N$ such that $q(x,x+\omega)$ is a strictly positive non-decreasing polynomial on $\N_0+N$ for all $\omega\in\Omega$. Moreover, $A_j=\{\om\in\cT_-\colon \om<j\}$ if $j\le0$, and $A_j=\{\om\in\cT_+\colon \om\ge j\}$ if $j>0$. It follows from Corollary \ref{co-5} that for all $x\in\N_0+N-\om_-$,

\begin{align}\nonumber
&T_{\pi}(x)\Bigl(\sum_{\omega\in A_0}q(x,x+\omega)+\sum_{ \omega\in A_1}q(x-1,x-1+\omega)\Bigr)\\ \nonumber
&\qquad\quad+\sum_{j=\omega_-}^{-1}T_{\pi}(x-j)\cdot\Bigl(\sum_{\om\in A_{j}}
q(x-j,x-j+\omega)-\sum_{\om\in A_{j+1}}
q(x-(j+1),x-(j+1)+\omega)\Bigr)\\ \label{Eq-1}
&= \sum_{j=1}^{\omega_+}
T_{\pi}(x-j)\Bigl(\sum_{\om\in A_j}
q(x-j,x-j+\omega)-\sum_{\om \in A_{j+1}}
q(x-(j+1),x-(j+1)+\omega)\Bigr).
\end{align}
Furthermore, {note that $R_->R_+$, and we have the following estimates for both sides of the above equality:
\[\begin{split}
  \text{LHS}_{\eqref{Eq-1}}&=T_{\pi}(x)\Bigl(\sum_{\omega\in A_0}q(x,x+\omega)+\sum_{\omega\in A_1}q(x-1,x-1+\omega)\Bigr)
  +\sum_{j=\omega_-}^{-1}T_{\pi}(x-j)\\
&\quad\cdot\Bigl(-q(x-(j+1),x-(j+1)+{j})\\&+\sum_{\om\in A_{j}}\Bigl(q(x-j,x-j+\omega)-q(x-(j+1),x-(j+1)+\omega)\Bigr)\Bigr)\\
&  \le T_{\pi}(x)\Bigl(\sum_{\omega\in A_0}q(x,x+\omega)+
\sum_{\omega\in A_1}q(x-1,x-1+\omega)\\ &+\sum_{j=\omega_-}^{-1}\sum_{\om\in A_{j}}\bigl(q(x-j,x-j+\omega)-q(x-(j+1),x-(j+1)+\omega)\bigr)\Big)\\
&=T_{\pi}(x)x^{R_-}\Big(\alpha_{0}+\rO\Big(x^{-\tilde{\sigma}}\Big)\Big),
\end{split}\]}where $\tilde{\sigma}=\min\{1,R_--R_+\}>0$.

{By the monotonicity of $q(x,x+\omega)$,
\[\begin{split}
\text{RHS}_\eqref{Eq-1}&\ge\sum_{j=1}^{\omega_+}
T_{\pi}(x-j)\Bigl(\sum_{\om\in A_j}
q(x-j,x-j+\omega)-\sum_{\om \in A_{j+1}}
q(x-j,x-j+\omega)\Bigr)\\
&= T_{\pi}(x-1)\Bigl(\sum_{\om\in A_1}q(x,x+\omega)-\sum_{\om\in A_1}\bigl(q(x,x+\omega)-q(x-\om,x)\bigr)\Bigr)\\
&\ge T_{\pi}(x-1)\lt(\alpha_1x^{R_+}+\rO\lt(x^{R_+-1}\rt)\rt).\end{split}\]}
Then there exist $N_1>N_2>0$ such that for all $x\ge N_1$,
\eqb\label{Eq-7}
\frac{T_{\pi}(x)}{T_{\pi}(x-1)}\ge x^{R_+-R_-}\frac{\alpha_1+\rO(x^{-1})}{\alpha_{0}+\rO(x^{-\tilde{\sigma}})}
=x^{R_+-R_-}\lt(\frac{\alpha_1}{\alpha_{0}}+\rO(x^{-\tilde{\sigma}})\rt)\ge \frac{\alpha_1}{\alpha_{0}}x^{R_+-R_-}\lt(1-N_2x^{-\tilde{\sigma}}\rt).
\eqe
Hence if $0<R_--R_+<1$, then $\tilde{\sigma}=R_--R_+<1$, and there exists $\widetilde{C}=\widetilde{C}(N_1,N_2)>0$ such that for all $x\ge N_1$,
\begin{align*}
T_{\pi}(x)&\ge T_{\pi}(N_1-1)\prod_{j=0}^{x-N_1}\lt(\frac{\alpha_1}{\alpha_{0}}(x-j)^{R_+-R_-}\lt(1-\frac{N_2}{(x-j)^{\tilde{\sigma}}}\rt)\rt)\\
&=T_{\pi}(N_1-1)\prod_{j=N_1}^{x}\lt(\frac{\alpha_1}{\alpha_{0}}j^{R_+-R_-}\rt)\prod_{j=N_1}^{x}\lt(1-\frac{N_2}{j^{\tilde{\sigma}}}\rt)\\
&\gtrsim \lt(\frac{\alpha_1}{\alpha_{0}}\rt)^{x+1-N_1}\frac{\Gamma\lt(x+1\rt)^{R_+-R_-}}{\Gamma(N_1)^{R_+-R_-}}\exp\Big(\sum_{j=N_1}^x-2N_2
j^{-\tilde{\sigma}}\Big) \\
&\gtrsim\Gamma(x)^{R_+-R_-}\lt(\frac{\alpha_1}{\alpha_{0}}\rt)^{x-\widetilde{C}x^{1-\tilde{\sigma}}}x^{R_+-R_-},\end{align*}
since $1-x^{-1}\ge \exp(-2x^{-1})$ for large $x$, and
we employ the fact that $T_{\pi}(N_1-1)>0$ since $\supp\pi=\cX$ is unbounded. Hence \eqref{Eq-27-a} holds.
Similarly, if $R_--R_+\ge1$, then {$\widetilde{\sigma}=1$}, and analogous arguments can be applied.

Next we show \eqref{Eq-27-b}.
Rewrite \eqref{Eq-46-a},
 \eqb\label{Eq-46-c}\sum_{j=\omega_-+1}^0\lt(\alpha_{j}+\beta_{j}(x)\rt)\pi\lt(x-j\rt)=x^{R_+-R_-}
\sum_{j=1}^{\omega_+}\lt(\alpha_j+\beta_j(x)\rt)\pi\lt(x-j\rt),\eqe
Summing up in \eqref{Eq-46-c} from $x$ to infinity yields \eqb\label{Eq-4-a}
\sum_{y=x}^{\infty}\sum_{j=\omega_-+1}^0\lt(\alpha_{j}+\beta_{j}(y)\rt)\pi\lt(y-j\rt)
=\sum_{y=x}^{\infty}y^{R_+-R_-}
\sum_{j=1}^{\omega_+}\lt(\alpha_j+\beta_j(y)\rt)\pi\lt(y-j\rt).\eqe
{Similarly, we can obtain the other estimates. It follows from \eqref{Eq-38} that there exists $C=C(N_4)>0$ and $N_5\in\N$ such that for all $x\ge N_5$,
\[\begin{split}
\text{LHS}_{\eqref{Eq-4-a}}
&\ge\lt(\alpha_0-Cx^{-\sigman}\rt)T_{\pi}(x),
\end{split}\]
\[\begin{split}
\text{RHS}_{\eqref{Eq-4-a}}&\le \sum_{y=x}^{\infty}y^{R_+-R_-}\sum_{j=1}^{\omega_+}\lt(\alpha_j+N_4y^{-\sigman}\rt)\pi\lt(y-j\rt)\\
& \le x^{R_+-R_-}T_{\pi}(x-\omega_+)\lt(\alpha_++Cx^{-\sigman}\rt),\end{split}\]}
which together further imply that
\[\frac{T_{\pi}(x)}{T_{\pi}(x-\omega_+)}\le x^{R_+-R_-}\frac{\alpha_++Cx^{-\sigman}}{\alpha_{0}-Cx^{-\sigman}}=x^{R_+-R_-}\lt(\frac{\alpha_+}{\alpha_{0}}+
\rO(x^{-\sigman})\rt).\]
The remaining arguments are analogous to the arguments for \eqref{Eq-27-a}.

\medskip

{\rm(ii)} $R_-=R_+$ and $\alpha_->\alpha_+$.
Analogous to (i), we will show that there exist real constants $\delta_+,\delta_-$ and $\widetilde{C}>0$ such that for all $\overline{\delta}>\delta_+$ and $\underline{\delta}<\delta_-$,
\begin{align}
 T_{\pi}(x)&\gtrsim\lt(\frac{\alpha_+}{\alpha_-}\rt)^{x+\widetilde{C}x^{1-\sigman}+\rO(\log x)}, \label{SEq-31-a} \\
 T_{\pi}(x)&\lesssim\lt(\frac{\alpha_+}{\alpha_-}\rt)^{(\omega_+-\omega_--1)^{-1}x+\widetilde{C}x^{1-\sigman}+\rO(\log x)}.
\label{SEq-31-b}
\end{align}
\medskip

We first prove \eqref{SEq-31-a}.
Since $R=R_-=R_+$, $$f_j(x)=x^R(\alpha_j+\beta_j(x)),\quad \beta_j(x)=\rO(x^{-\sigman}),\quad  j=\omega_-+1,\ldots,\om_+.$$
 Moreover, $\alpha=\alpha_+-\alpha_-<0$ implies that {$$\sum_{j=1}^{\omega_+}\alpha_j<\sum_{j=\omega_-+1}^0\alpha_{j}$$}
 From {\eqref{Eq-46-a}}, it follows that \[\sum_{j=\omega_-+1}^0\pi\lt(x-j\rt)\alpha_{j}+\sum_{j=\omega_-+1}^0\pi\lt(x-j\rt)\beta_{j}(x)
=\sum_{j=1}^{\omega_+}\pi\lt(x-j\rt)\alpha_{j}+\sum_{j=1}^{\omega_+}\pi\lt(x-j\rt)\beta_{j}(x).
\]
Summing up the above equality from $x$ to $\infty$ yields
\[\begin{split}
&\sum_{j=\omega_-+1}^0\sum_{y=x}^{\infty}\pi\lt(y-j\rt)\alpha_{j}+\sum_{j=\omega_-+1}^0\sum_{y=x}^{\infty}\pi\lt(y-j\rt)\beta_{j}(y)\\
&\quad=\sum_{j=1}^{\omega_+}\sum_{y=x}^{\infty}\pi\lt(y-j\rt)\alpha_j+\sum_{j=1}^{\omega_+}\sum_{y=x}^{\infty}\pi\lt(y-j\rt)\beta_j(y).
\end{split}\]
Since each double sum in the above equality is convergent, we have \begin{align*}
  0&=\sum_{j=\omega_-+1}^0\alpha_{j}\sum_{y=x}^{\infty}(\pi(y)-\pi\lt(y-j\rt))+\sum_{j=1}^{\omega_+}\alpha_j\sum_{y=x}^{\infty}(\pi\lt(y-j\rt)-\pi(y))\\
  &\quad -\lt(\sum_{j=\omega_-+1}^0\alpha_{j}-\sum_{j=1}^{\omega_+}\alpha_j\rt)T_{\pi}(x)\\
  &\quad+\sum_{j=\omega_-+1}^0\sum_{y=x}^{\infty}(\pi(y)\beta_{j}(y+j)-\pi\lt(y-j\rt)\beta_{j}(y))\\
  &\quad+\sum_{j=1}^{\omega_+}\sum_{y=x}^{\infty}(\pi\lt(y-j\rt)\beta_j(y)-\pi(y)\beta_{j}(y+j))\\
  &\quad+\sum_{j=1}^{\omega_+}\sum_{y=x}^{\infty}\pi(y)\beta_j(y+j)-\sum_{j=\omega_-+1}^0\sum_{y=x}^{\infty}\pi(y)\beta_{j}(y+j).
\end{align*}
This further yields the following equality
\begin{align}\nonumber
 &(\alpha_--\alpha_+)T_{\pi}(x)+\sum_{j=\omega_-+1}^0\sum_{y=x}^{\infty}\pi(y)\beta_{j}(y+j)-
\sum_{j=1}^{\omega_+}\sum_{y=x}^{\infty}\pi(y)\beta_j(y+j)\\\label{SEq-43}
&=\sum_{j=\omega_-+1}^{-1}\sum_{\ell=0}^{-j-1}\pi\lt(x+\ell\rt)\widetilde{f}_{j}(x+j+\ell)+\sum_{j=1}^{\omega_+}
 \sum_{\ell=1}^{j}\pi\lt(x-\ell\rt)\widetilde{f}_{j}(x+j-\ell),
\end{align}
where $\widetilde{f}_j(x)=\alpha_j+\beta_j(x)=x^{-R}f_j(x)\ge0$, $j=\om_-+1,\ldots,\om_+$.
From \eqref{Eq-38}, it follows that there exist $C,N>0$ such that
\[\Bigl|\sum_{j=\omega_-+1}^0\sum_{y=x}^{\infty}\pi(y)\beta_{j}(y+j)-
\sum_{j=1}^{\omega_+}\sum_{y=x}^{\infty}\pi(y)\beta_j(y+j)\Bigl| \
\le C\sum_{y=x}^{\infty}\pi(y)y^{-\sigman}\le Cx^{-\sigman}T_{\pi}(x),\]
for $x\ge N$. Hence
\[\text{LHS}_{\eqref{SEq-43}}\le\lt((\alpha_--\alpha_+)+Cx^{-\sigman}\rt)T_{\pi}(x).\]
Using Fubini's theorem, we have:
\eqb\label{Eq-46}
\text{RHS}_{\eqref{SEq-43}}=\sum_{j=\om_-+2}^0\pi(x-j)\sum_{\ell=\om_-+1}^{j-1}\widetilde{f}_{\ell}(x+\ell-j)+
\sum_{j=1}^{\om_+}\pi(x-j)\sum_{\ell=j}^{\om_+}\widetilde{f}_{\ell}(x+\ell-j).
\eqe
Hence further choosing larger $N$ and $C$, we have for all $x\ge N$,

\begin{align}\nonumber
\text{RHS}_{\eqref{SEq-43}} &\ge\pi\lt(x-1\rt)\sum_{j=1}^{\omega_+}\widetilde{f}_j(x+j-1)
 =\pi\lt(x-1\rt)\sum_{j=1}^{\omega_+}(\alpha_{j}+\beta_{j}(x+j-1))\\
 &\ge (\alpha_+-Cx^{-\sigman})\pi\lt(x-1\rt), \label{Eq-45}
\end{align}
This implies from $\pi(x-1)=T_{\pi}(x-1)-T_{\pi}(x)$ that  
\[\frac{T_{\pi}(x)}{T_{\pi}(x-1)}\ge\frac{\alpha_+-Cx^{-\sigman}}{\alpha_-}.\] 
Using similar arguments as in the proof of (i), one obtains \eqref{SEq-31-a}.

Next we show \eqref{SEq-31-b}.
We establish the reverse estimates for both sides of \eqref{SEq-43}.
Similarly, there exists some $C,N>0$ such that for all $x\ge N$,
\[
  \text{LHS}_{\eqref{SEq-43}}\ge(\alpha_--\alpha_+-Cx^{-\sigman})T_{\pi}(x)\ge(\alpha_--\alpha_+-Cx^{-\sigman})T_{\pi}(x-\om_--1),\]
\[
\text{RHS}_{\eqref{SEq-43}}\le(\alpha_++Cx^{-\sigman})\lt(T_{\pi}(x-\omega_+)-T_{\pi}(x-\omega_--1)\rt).\]
 This implies that for a possibly larger $C,\ N$, for all $x\ge N$,
\[\frac{T_{\pi}(x-\omega_--1)}{T_{\pi}(x-\omega_+)}\le\frac{\alpha_++Cx^{-\sigman}}{\alpha_-}.\]
The remaining arguments are similar to those in the proof of (i).

\smallskip

{\rm(iii)-(vi)} $\alpha=0$. Hence $\alpha_+=\alpha_-$. Recall $$\Delta=\alpha_+^{-1}\cdot\cab-\gamma,\qquad\quad\ \, \text{if}\ \sigman<1,\\ -\gamma+R\vartheta,\quad \text{if}\ \sigman=1.\cae$$
Let $\delta=\Delta(\omega_+-\omega_--1)^{-1}$. 
For $j=\om_-+1,\ldots,\om_+$, $$r_j=\cab \gamma_j,\qquad\qquad\ \, \text{if}\ \sigman<1,  \\ \gamma_j-jR\alpha_j,\quad \text{if}\ \sigman=1,\cae$$
$$\vartheta_j(x)=\beta_j(x)-r_j x^{-\sigman}.$$
Hence we have
$$\vartheta_j(x)=\rO(x^{-\overline{\rhon}}),\quad j=\omega_-+1,\ldots,\om_+,$$
where
$$\overline{\rhon}=\left\{\begin{array}{ll}
\min\{1,\rhon\},& \text{if}\ \sigman<1,\\
 \rhon,& \text{if}\ \sigman=1.\end{array}\right.$$
 Let $\eta=\overline{\rhon}-\sigman$ and $\varepsilon=\min\{\sigman,\eta\}$. Hence $0<\varepsilon\le\eta\le1$. If $\sigman<1$, then $\eta\le1-\sigman$.  If $\sigman=1$, then $\varepsilon=\eta$.

To show (iii)-(vi), it suffices to prove that there exists $C>0$ such that
\begin{equation}\label{Eq-6-a}
T_{\pi}(x)\gtrsim\left\{\begin{array}{ll}
\exp\Bigl(-\frac{\Delta}{1-\sigman} x^{1-\sigman}+\rO\bigl(x^{1-\sigman-\varepsilon}+\log x\bigr)\Bigr),& \text{if}\ \sigman<1,\ \sigman+\varepsilon\neq1, \Delta>0\\
   \exp\Bigl(-\frac{\Delta}{1-\sigman} x^{1-\sigman}+\rO(\log x)\Bigr),& \text{if}\ \sigman<1,\ \sigman+\varepsilon=1, \Delta>0\\
   x^{-(\Delta-1)},& \text{if}\ \sigman=1, \Delta>0,\\
 \exp\Bigl(-\frac{C}{1-\rhon}x^{1-\rhon}+\rO\bigl(x^\sigma+\log x\bigr)\Bigr), & \text{if}\ {\rhon}<1,\ \sigman+\rhon\neq1, \Delta=0\\
   \exp\Bigl(-\frac{C}{1-\rhon}x^{1-{\rhon}}+\rO(\log x)\Bigr), & \text{if}\ \rhon<1,\ \sigman+\rhon=1, \Delta=0\\
   x^{-C},& \text{if}\ {\rhon}\ge1, \Delta=0,
 \end{array}\right.
\end{equation}
where $\sigma=\max\{1-\sigman-\rhon,
   0\}$,  and if $\Delta>0$, then
\begin{equation}\label{Eq-6-b}
T_{\pi}(x)\lesssim \left\{\begin{array}{ll} \exp\Bigl(-\frac{\delta}{1-\sigman} x^{1-\sigman}+\rO\bigl(x^{\max\{1-\sigman-\varepsilon,0\}}+\log x\bigr)\Bigr), &  \text{if}\ \sigman<1,\ \sigman+\varepsilon\neq1,\\
   \exp\Bigl(-\frac{\delta}{1-\sigman} x^{1-\sigman}+\rO(\log x)\Bigr),& \text{if}\ \sigman<1,\ \sigman+\varepsilon=1,\\
   x^{-\max\{\delta-1,0\}},& \text{if}\ \sigman=1.
   \end{array}\right.
   \end{equation}

To show \eqref{Eq-6-a} and \eqref{Eq-6-b} for a probability distribution $\mu$ on $\N_0$, define its {\em weighted} tail distribution on $\N_0$ as
$$W_{\mu}\colon\N\to[0,1],\quad x\mapsto\sum_{y=x}^{\infty}y^{-\sigman}\mu(y).$$
In the following, we will show there exist constants $C>\sigman$ such that
\begin{equation}\label{SEq-33-a}
W_{\pi}(x)\gtrsim \left\{\begin{array}{ll}
\exp\Bigl(\frac{-\Delta}{1-\sigman}x^{1-\sigman}+\rO\bigl(x^{1-\sigman-\varepsilon}\bigr)\Bigr),& \text{if}\ \sigman<1,\ \sigman+\varepsilon\neq1,  \Delta>0,\\
   \exp\Bigl(\frac{-\Delta}{1-\sigman}x^{1-\sigman}+\rO(\log x)\Bigr),& \text{if}\ \sigman<1,\ \sigman+\varepsilon=1,  \Delta>0, \\
   x^{-\Delta},& \text{if}\ \sigman=1,  \Delta>0,\\
   \exp\Bigl(\frac{-C}{1-\rhon}x^{1-\rhon}+\rO\bigl(x^{\max\{1-\sigman-{\rhon},
   0\}}\bigr)\Bigr),& \text{if}\ {\rhon}<1,\ \sigman+\rhon\neq1,  \Delta=0,\\
   \exp\Bigl(\frac{-C}{1-{\rhon}}x^{1-{\rhon}}+\rO(\log x)\Bigr),& \text{if}\ \rhon<1,\ \sigman+\rhon=1,  \Delta=0,\\
   x^{-C},& \text{if}\ {\rhon}\ge1, \Delta=0,
 \end{array}\right.
 \end{equation}
with $\Delta>1$, when $\sigman=1$. Moreover, if $\Delta>0$, then
\begin{equation}\label{SEq-33-b}
W_{\pi}(x)\lesssim \left\{\begin{array}{ll}
\exp\Bigl(\frac{-\delta}{1-\sigman}x^{1-\sigman}+\rO(x^{\max\{1-\sigman-\varepsilon,0\}})\Bigr),& \text{if}\ \sigman<1,\ \sigman+\varepsilon\neq1,\\
   \exp\Bigl(\frac{-\delta}{1-\sigman}x^{1-\sigman}+\rO(\log x)\Bigr),& \text{if}\ \sigman<1,\ \sigman+\varepsilon=1,\\
   x^{-\delta},& \text{if}\ \sigman=1.
   \end{array}\right.
   \end{equation}
Then we will prove \eqref{Eq-6-a} and \eqref{Eq-6-b} based on \eqref{SEq-33-a} and \eqref{SEq-33-b}.

\medskip

\noindent\textbf{Step I.} Prove \eqref{SEq-33-a} and \eqref{SEq-33-b}. Here we will also show $\Delta\ge0$, and in particular $\Delta>1$ when $\sigman=1$.
We first show  \eqref{SEq-33-a}.
Since $\alpha_+=\alpha_-$,
\[\begin{split}
\text{LHS}_{\eqref{SEq-43}}&=\sum_{j=\omega_-+1}^0\sum_{y=x}^{\infty}\pi(y)\lt(r_j(y+j)^{-\sigman}+\vartheta_{j}(y+j)\rt)\\&\quad-\sum_{j=1}^
{\omega_+}\sum_{y=x}^{\infty}\pi(y)\lt(r_j(y+j)^{-\sigman}+\vartheta_{j}(y+j)\rt)\\
&=\lt(\sum_{j=\omega_-+1}^0r_j-\sum_{j=1}^{\omega_+}r_j\rt)W_{\pi}(x)\\&\quad+\sum_{j=\omega_-+1}^0\sum_{y=x}^{\infty}\pi(y)\lt(\vartheta_{j}(y+j)+r_j((y+j)^
{-\sigman}-y^{-\sigman})\rt)\\
&\quad-\sum_{j=1}^{\omega_+}\sum_{y=x}^{\infty}\pi(y)\lt(\vartheta_{j}(y+j)+r_j((y+j)^{-\sigman}-y^{-\sigman})\rt).
\end{split}\]
By Lemma~\ref{Sle-1},
\[\sum_{j=\omega_-+1}^0r_j-\sum_{j=1}^{\omega_+}r_j=\alpha_+\Delta.\]
Moreover,
\[\begin{split}
&\hspace{-2cm}\Bigl|\sum_{j=\omega_-+1}^0\sum_{y=x}^{\infty}\pi(y)\lt(\vartheta_{j}(y+j)+r_j((y+j)^{-\sigman}-y^{-\sigman})\rt)\\
&\quad-\sum_{j=1}^{\omega_+}\sum_{y=x}^{\infty}\pi(y)\lt(\vartheta_{j}(y+j)+r_j((y+j)^{-\sigman}-y^{-\sigman})\rt)\Bigl|\\
 & \lesssim \sum_{y=x}^{\infty}\pi(y)y^{-\sigman-\eta}\le x^{-\eta}W_{\pi}(x).
\end{split}\]
Since $\text{RHS}_{\eqref{SEq-43}}\ge0$ for all large $x$, we have $\Delta\ge0$.

From \eqref{Eq-45} it follows that there exist $C,\ N\in\N$ such that for all $x\ge N$,
\[
  \text{LHS}_{\eqref{SEq-43}}\le\alpha_+\lt(\Delta+Cx^{-\eta}\rt)W_{\pi}(x),\]
while\[
\begin{split}
\text{RHS}_{\eqref{SEq-43}}&\ge\alpha_+(1-Cx^{-\sigman})\pi(x-1)\\
&=\alpha_+(1-Cx^{-\sigman})(x-1)^{\sigman}(W_{\pi}(x-1)-W_{\pi}(x))\\
&=\alpha_+(x^{\sigman}-C-\sigman x^{\sigman-1}+\rO(x^{\sigman-2}))(W_{\pi}(x-1)-W_{\pi}(x)).
\end{split}\]
Further choosing larger $N$ and $C$, by the monotonicity of $W_{\pi}$, for all $x\ge N$,
\begin{multline*}
(x^{\sigman}-C-\sigman x^{\sigman-1}+\rO(x^{\sigman-2}))W_{\pi}(x-1)\\
\le \lt(x^{\sigman}-C+\Delta-\sigman x^{\sigman-1}+Cx^{-\eta}+\rO(x^{\sigman-2})\rt)W_{\pi}(x).\end{multline*}
If $\sigman<1$, then $\eta\le1-\sigman$, and hence $\eta+\sigman-2\le-\sigman$. If $\sigman=1$, then $\varepsilon=\eta$. Recall $\overline{\rhon}=\eta+\sigman$. Then we have
  \[\begin{split}
  \frac{W_{\pi}(x)}{W_{\pi}(x-1)}&\ge\frac{x^{\sigman}-C-\sigman x^{\sigman-1}+\rO(x^{\sigman-2})}
  {x^{\sigman}-C+\Delta+Cx^{-\eta}-\sigman  x^{\sigman-1}+\rO(x^{\sigman-2})}\\
  &=\left\{\begin{array}{ll}
  1-\Delta x^{-\sigman}(1+\rO(x^{-\varepsilon})),& \text{if}\ \Delta>0,\\ 1-Cx^{-\overline{\rhon}}(1+\rO(x^{\max\{-\sigman,\overline{\rhon}-2\}})),& \text{if}\ \Delta=0.\end{array}\right.
  \end{split}\]
 First assume $\Delta>0$. Since $\varepsilon\le1$, by Euler-Maclaurin's formula,
  \[\begin{split}
  \log\frac{W_{\pi}(x)}{W_{\pi}(N-1)}&\ge\sum_{j=N}^{x}\log(1-\Delta j^{-\sigman}+\rO(j^{-\sigman-\varepsilon}))\\
  &=\left\{\begin{array}{ll}\frac{-\Delta}{1-\sigman}x^{1-\sigman}+\rO(x^{\max\{0,1-\sigman-\varepsilon\}}), & \text{if}\ \sigman<1,\ \sigman+\varepsilon\neq1,\\
   \frac{-\Delta}{1-\sigman}x^{1-\sigman}+\rO(\log x),& \text{if}\ \sigman<1,\ \sigman+\varepsilon=1,\\
   -\Delta\log x+\rO(1),& \text{if}\ \sigman=1,\end{array}\right.
  \end{split}\]
 which implies that
  \[
   W_{\pi}(x)\gtrsim\left\{\begin{array}{ll}\exp\Bigl(\frac{-\Delta}{1-\sigman}x^{1-\sigman}+\rO(x^{1-\sigman-\varepsilon})\Bigr),& \text{if}\ \sigman<1,\ \sigman+\varepsilon\neq1,\\
   \exp\Bigl(\frac{-\Delta}{1-\sigman}x^{1-\sigman}+\rO(\log x)\Bigr),& \text{if}\ \sigman<1,\ \sigman+\varepsilon=1,\\
   x^{-\Delta},& \text{if}\ \sigman=1,\end{array}\right.
 \]
  i.e., \eqref{SEq-33-a} holds. Moreover, since $x^\sigman W_{\pi}(x)\le T_{\pi}(x)\to0$ as $x\to\infty$, we have
  $$\Delta>\left\{\begin{array}{ll}0,& \text{if}\ \sigman<1,\\ 1,& \text{if}\ \sigman=1.\end{array}\right.$$

  Now assume $\Delta=0$, then
  \[
   W_{\pi}(x)\gtrsim\left\{\begin{array}{ll}\exp\Bigl(\frac{-C}{1-\overline{\rhon}}x^{1-\overline{\rhon}}+\rO(x^{\max\{1-\sigman-\rhon,0\}})\Bigr),
  & \text{if}\ \overline{\rhon}\neq1,\ \sigman+\rhon\neq1,\\
   \exp\Bigl(\frac{-C}{1-\overline{\rhon}}x^{1-\overline{\rhon}}+\rO(\log x)\Bigr),&\text{if}\ \overline{\rhon}\neq1,\ \sigman+\rhon=1,\\
   x^{-C},& \text{if}\ \overline{\rhon}=1,\end{array}\right.
 \]
 where we use the fact that $\sigman+\overline{\rhon}=1$ implies $0<\sigman,\rhon<1$ and $\sigman+\rhon=1$. Moreover, also due to $x^\sigman W_{\pi}(x)\le T_{\pi}(x)\to0$ as $x\to\infty$, we have $\overline{\rhon}\le1$, which implies that $\sigman<1$. In addition, $C>\sigman$ when $\overline{\rhon}=1$, i.e., $\sigman<1\le\rhon$. Hence for some $C>\sigman$, \[
   W_{\pi}(x)\gtrsim\left\{\begin{array}{ll}\exp\Bigl(\frac{-C}{1-\rhon}x^{1-\rhon}+\rO(x^{\max\{1-\sigman-{\rhon},
   0\}})\Bigr),& \text{if}\ {\rhon}<1,\ \sigman+\rhon\neq1,\\
   \exp\Bigl(\frac{-C}{1-{\rhon}}x^{1-{\rhon}}+\rO(\log x)\Bigr),&  \text{if}\ \rhon<1,\ \sigman+\rhon=1,\\
   x^{-C},& \text{if}\ {\rhon}\ge1.\end{array}\right.
 \]

\smallskip

Next we show \eqref{SEq-33-b} by establishing the reverse estimates for both sides of \eqref{SEq-43}. From \eqref{Eq-46} it follows that there exist positive constants $N$ and $C_i$ ($i=1,2$) such that $x\ge N$,
\[
  \text{LHS}_{\eqref{SEq-43}}\ge\alpha_+\lt(\Delta-Cx^{-\eta}\rt)W_{\pi}(x)\ge\alpha_+\lt(\Delta-Cx^{-\eta}\rt)
  W_{\pi}(x-(\om_-+1)),\]
whereas
\[\begin{split}
\text{RHS}_{\eqref{SEq-43}}&\le\alpha_+\lt(1+C_1x^{-\sigman}\rt)\sum_{j=\om_-+2}^{\omega_+}\pi\lt(x-j\rt)\\
&\le\alpha_+\lt(x^{\sigman}+C_2\rt)\sum_{j=\om_-+2}^{\omega_+}\pi\lt(x-j\rt)(x-j)^{-\sigman}\\
&=\alpha_+\lt(x^{\sigman}+C_2\rt)(W_{\pi}(x-\omega_+)-W_{\pi}(x-(\omega_-+1))),
\end{split}\]
 Hence, when $\Delta>0$, then for all $x\ge N$, \[\frac{W_{\pi}(x-(\omega_-+1))}{W_{\pi}(x-\omega_+)}\le\frac{x^{\sigman}+C_2}{
  x^{\sigman}+C_2+\Delta-Cx^{-\eta}}=1-\Delta x^{-\sigman}+\rO(x^{-\sigman-\varepsilon}).\]
  Analogous to the above analysis, one might show
  \[
   W_{\pi}(x)\lesssim\left\{\begin{array}{ll}\exp\Bigl(\frac{-\delta}{1-\sigman}x^{1-\sigman}+\rO(x^{\max\{1-\sigman-\varepsilon,0\}})\Bigr),\quad \text{if}\ \sigman<1,\ \sigman+\varepsilon\neq1,\\
   \exp\Bigl(\frac{-\delta}{1-\sigman}x^{1-\sigman}+\rO(\log x)\Bigr),\qquad\qquad\quad\ \ \! \text{if}\ \sigman<1,\ \sigman+\varepsilon=1,\\
   x^{-\delta},\qquad\qquad\qquad\qquad\qquad\qquad\qquad\qquad \text{if}\ \sigman=1,\end{array}\right.
 \]
  in particular, one can further choose $\delta\ge1$ (not necessarily $\delta=\Delta(\omega_+-\omega_--1)^{-1}$) when $\sigman=1$, also due to $x^{\sigman} W_{\pi}(x)\le T_{\pi}(x)\to0$ as $x\to\infty$.

Moreover, one can always show $W_{\pi}(x)\le x^{-\sigman}T_{\pi}(x)\le x^{-\sigman}$, hence \eqref{SEq-33-b} also holds when $\sigman=1$.

\medskip

\noindent\textbf{Step II.} Prove \eqref{Eq-6-a} and \eqref{Eq-6-b} based on \eqref{SEq-33-a} and \eqref{SEq-33-b}.

Since $W_{\pi}(x)\le x^{-\sigman}T_{\pi}(x)$, \eqref{Eq-6-a} follows directly from  \eqref{SEq-33-a}.

Next, we prove \eqref{Eq-6-b} based on  \eqref{SEq-33-b}. Recall that $$\pi(x)\le x^{\sigman} W_{\pi}(x).$$ Assume $\Delta>0$. We only prove the case $\sigman<1$ and $\sigman+\varepsilon=1$. The other two cases can be proved using analogous arguments. Then there exist $N\in\N$ and $C_1>\sigman$ such that $\exp\Bigl(\frac{-\delta}{1-\sigman}y^{1-\sigman}+C_1\log y\Bigr)$ is decreasing on $[N,+\infty)$, and for all $x\ge N$,
{\begin{align*}
T_{\pi}(x)&=\sum_{y=x}^{\infty}\pi(y)\ \le \ \sum_{y=x}^{\infty}y^{\sigman} W_{\pi}(y)  \lesssim\sum_{y=x}^{\infty}\exp\Bigl(\tfrac{-\delta}{1-\sigman}y^{1-\sigman}+C_1\log y\Bigr)\\
&\lesssim\int_{x-1}^{\infty}\exp\Bigl(\tfrac{-\delta}{1-\sigman}y^{1-\sigman}+C_1\log y\Bigr){\rm d}y \lesssim\int_{x-1}^{\infty}y^{C_1+\sigman}{\rm d}\lt(-\exp\Bigl(\tfrac{-\delta}{1-\sigman}y^{1-\sigman}\Bigr)\rt)\\
&\le(x-1)^{C_1+\sigman}\exp\Bigl(\tfrac{-\delta}{1-\sigman}(x-1)^{1-\sigman}\Bigr)\\
&\quad+(C_1+\sigman)(x-1)^{\sigman-1}\int_{x-1}^{\infty}\exp\Bigl(\tfrac{-\delta}{1-\sigman}y^{1-\sigman}+C_1\log y\Bigr){\rm d}y,
\end{align*}}
which further implies that  for all $x\ge N$,
\[\begin{split}
T_{\pi}(x)&\lesssim\frac{(x-1)^{C_1+\sigman}\exp\Bigl(\frac{-\delta}{1-\sigman}(x-1)^{1-\sigman}\Bigr)}
{1+\rO\bigl((x-1)^{\sigman-1}\bigr)}\\
&\lesssim\exp\Bigl(\tfrac{-\delta}{1-\sigman}(x-1)^{1-\sigman}+\rO(\log x)\Bigr)=\exp\Bigl(\tfrac{-\delta}{1-\sigman}x^{1-\sigman}+\rO(\log x)\Bigr).\end{split}\]
This shows \eqref{Eq-6-b} in that case.

\section{Proof of Theorem~\ref{th-19b}}

Again ($\rm\mathbf{A3}$) implies $\supp\nu=\partial^{\sf c}$ \cite{CMM13}, which is unbounded by ($\rm\mathbf{A2}$).

Comparing the identities for stationary distributions and QSDs, the unique difference comes from an extra term on the RHS of the identity of QSDs with coefficient $\theta_{\nu}>0$. This makes the identity in Theorem~\ref{th-18}(3) for stationary distributions to be an inequality with its LHS greater than its RHS for QSDs. Hence, all arguments in the proof of Theorem~\ref{th-19b} establishing $\alpha\le0$ as well as the lower estimates for $T_{\pi}$ (the tail of the stationary distribution) carry over to $T_{\nu}$.

Next, we show $R\ge0$. The proof is in a similar spirit to that for $\alpha\le0$. Since $\alpha\le0$,  {then $R_-=R$ from the definition of $\alpha$.}  Again, assume w.o.l.g. that $\om_*=1$ such that $\partial^{\sf}$ contains all large positive integers by Proposition~\ref{pro-0}. From Theorem~\ref{th-18b}(3), similar to \eqref{Eq-46-a}, we have for all large $x$,
\begin{align*}
x^{R}(\alpha_-+Cx^{-\sigman})T_{\nu}(x)&\ge x^{R}\sum_{j=\omega_-+1}^0\lt(\alpha_{j}+\beta_{j}(x)\rt)\nu\lt(x-j\rt)\\
&= \theta_{\nu}T_{\nu}(x)+x^{R_+}\sum_{j=1}^{\omega_+}\lt(\alpha_j+\beta_j(x)\rt)\nu\lt(x-j\rt) \\
&\ge\theta_{\nu}T_{\nu}(x),\end{align*}
which yields\[x^{R}(\alpha_-+Cx^{-\sigman})-\theta_{\nu}\ge0,\]
This shows $R\ge0$, since $\theta_{\nu}>0$. Moreover, if $R=0$, then $\alpha_-\ge\theta_{\nu}$. The claim that $R_-=R_+=0$ implies $\alpha\le-\theta_{\nu}$ is proved below in (vii).

{Recall that $\alpha_0=\beta$.} Similar to \eqref{Eq-1} and the inequality \eqref{Eq-7} based on it, one can also obtain $\alpha_0\ge\theta_{\nu}$ if $R=0$ and $R_->R_+$. Moreover, there exists $C>\alpha_1>0$ such that for all large $x$,
\eqb\label{Eq-7-a}\frac{T_{\nu}(x)}{T_{\nu}(x-1)}\ge \left\{\begin{array}{ll} x^{R_+-R_-}\frac{\alpha_1-Cx^{-1}}{\alpha_0-\theta_{\nu}x^{-R_-}+Cx^{-\tilde{\sigma}}},&\text{if}\ R_->R_+,\\
\frac{\alpha_1-Cx^{-1}}{\alpha_0+\alpha_1-\theta_{\nu}x^{-R}+Cx^{R-1}},& \text{if}\ R_-=R_+,\end{array}\right.\eqe
where we recall $\tilde{\sigma}=\min\{1,R_--R_+\}$.

Similar to \eqref{Eq-4-a}, we establish  \eqb\label{Eq-8}
\sum_{y=x}^{\infty}\sum_{j=\omega_-+1}^0\lt(\alpha_{j}+\beta_{j}(y)\rt)\nu\lt(y-j\rt)
=\theta_{\nu}\sum_{y=x}^{\infty}y^{-R_-}T_{\nu}(y)+\sum_{y=x}^{\infty}y^{R_+-R_-}
\sum_{j=1}^{\omega_+}\lt(\alpha_j+\beta_j(y)\rt)\nu\lt(y-j\rt).
\eqe
Since LHS$_{\eqref{Eq-8}}$ is finite, we have $\sum_{y=x}^{\infty}y^{-R_-}T_{\nu}(y)$ is also finite. Furthermore,  by similar analysis as in the proof of Theorem~\ref{th-19} that there exists $C>0$  such that for all large $x\in\N$,
\[\theta_{\nu}\sum_{y=x}^{\infty}y^{-R_-}T_{\nu}(y)\le (\alpha_-+Cx^{-\sigman})T_{\nu}(x)\to0,\quad \text{if}\ x\to\infty.\]

\noindent{\bf Step I}. Establish lower estimates for $T_{\nu}$ based on the above inequality, using similar asymptotic analysis demonstrated repeatedly in the proof of  Theorem~\ref{th-19b}.

\smallskip

\noindent{\rm (1)} $R_-=R>R_+$.

\smallskip
\noindent{$\bullet$} $R=0$ (cases (i)-(iii)).  Then $\alpha_{0}\ge\theta_{\nu}$. If $\alpha_{0}>\theta_{\nu}$, then there exists $\widetilde{C}>0$ such that \[ T_{\nu}(x)\gtrsim\exp\lt(-(R_--R_+)\log\Gamma(x)-\lt(\log\tfrac{\alpha_{0}-\theta_{\nu}}{\alpha_1}\rt)x-\widetilde{C}x^{1-(R_--R_+)}
+\rO(\log x)\rt),\]i.e., $\nu\in\cP^{1-}_{R_--R_+}$. Hence case (i) is proved.
If $\alpha_{0}=\theta_{\nu}$, then
\[\frac{T_{\nu}(x)}{T_{\nu}(x-1)}\ge x^{\min\{0,1+R_+-R_-\}}(\tfrac{\alpha_1}{C}-x^{-1}),\]
which yields that\[T_{\nu}(x)\gtrsim\exp\lt(\min\{0,R_+-R_-+1\}\log\Gamma(x)-\lt(\log\tfrac{C}{\alpha_1}\rt)x-\tfrac{C}{\alpha_1}\log x\rt),\]
i.e., $\nu\in\cP^{2-}_{1}$ if $0>R_+-R_-\ge-1$, and $\nu\in\cP^{1-}_{R_--R_+-1}$ if $R_+-R_-<-1$. Hence the cases (ii) and (iii) are also proved.

\smallskip
\noindent{$\bullet$} $R>0$ (case (iv)). Based on \eqref{Eq-7-a}, there exists $\widetilde{C}>0$ such that
\[\log T_{\nu}(x)\gtrsim \exp\lt(-(R_--R_+)\log\Gamma(x)-\lt(\log\tfrac{\alpha_{0}}{\alpha_1}\rt)x-\widetilde{C}x^{1-\min\{\widetilde{\sigma},R_-\}}
+\rO(\log x)\rt),\]
i.e., $\nu\in\cP^{1-}_{R_--R_+}$. Hence the former part of (iv) is proved. The second part of (iv) is proved below in Step II.
\smallskip

\noindent{\rm(2)} $R_-=R_+=R$. Then \eqref{Eq-8} is \[\sum_{y=x}^{\infty}\sum_{j=\omega_-+1}^0\lt(\alpha_{j}+\beta_{j}(y)\rt)\nu\lt(y-j\rt)
=\theta_{\nu}\sum_{y=x}^{\infty}y^{-R_-}T_{\nu}(y)+\sum_{y=x}^{\infty}
\sum_{j=1}^{\omega_+}\lt(\alpha_j+\beta_j(y)\rt)\nu\lt(y-j\rt).
\]
from which it implies that there exists $C>0$ and $N\in\N$ such that for all large $x\in\N$,
\eqb\label{Eq-5}
\frac{T_{\nu}(x)}{T_{\nu}(x-1)}\ge\frac{\alpha_+-Cx^{-\sigman}}{\alpha_--\theta_{\nu}x^{-R_-}+Cx^{-\sigman}},
\eqe

Based on which we establish the following lower estimates for $T_{\nu}(x)$.

\smallskip

\noindent{(v)} $R>0$ and $\alpha<0$.
We can show
$$T_{\nu}(x)\gtrsim\cab\exp\bigl((\log\frac{\alpha_+}{\alpha_-})x+\rO(\log x)\bigr),\qquad\quad\quad\ \, \text{if}\ \min\{R,\sigman\}=1,\\ \exp\bigl((\log\frac{\alpha_+}{\alpha_-})x+\rO(x^{1-\min\{R,\sigman\}})\bigr),\quad \text{if}\ \min\{R,\sigman\}\neq1,\cae$$i.e.,
$\nu\in\cP^{2-}_{1}$. The latter part is proved in Step II below.

\smallskip

\noindent{(vi)} $R>0$ and $\alpha=0$.
We prove the conclusions case by case.

\smallskip
\noindent{$\bullet$} $0<R<\sigman$ and $\alpha=0$. Then
$$1\ge T_{\nu}(x)\gtrsim\left\{\begin{array}{ll} \exp\bigl(\frac{\theta_{\nu}}{\alpha_+(1-R)}x^{1-R}+\rO(\log x)\bigr),& \text{if}\ \min\{2R,\sigman\}=1,\\
\exp\bigl(\frac{\theta_{\nu}}{\alpha_+(1-R)}x^{1-R}+\rO(x^{1-\min\{2R,\sigman\}})\bigr),& \text{if}\ \min\{2R,\sigman\}\neq1,\end{array}\right.$$
 which tends to infinity as $x\to\infty$.
This is a contradiction, and thus this case is not possible to occur.

\noindent{$\bullet$} $0<R=\sigman<1$ and $\alpha=0$. Then
$$T_{\nu}(x)\gtrsim\left\{\begin{array}{ll} \exp\bigl(-\frac{2C-\theta_{\nu}}{\alpha_+(1-R)}x^{1-R}+\rO(\log x)\bigr),& \text{if}\ 2R=1,\\  \exp\bigl(-\frac{2C-\theta_{\nu}}{\alpha_+(1-R)}x^{1-R}+\rO(x^{1-2R})\bigr),& \text{if}\ 2R\neq1,\end{array}\right.$$
i.e., $\nu\in\cP^{2-}_{1-R}$.

\noindent{$\bullet$} $\min\{1,R\}>\sigman$ and $\alpha=0$. Then $$T_{\nu}(x)\gtrsim\left\{\begin{array}{ll} \exp\bigl(-\frac{2C}{\alpha_+(1-\sigman)}x^{1-\sigman}+\rO(\log x)\bigr),& \text{if}\ \min\{R,2\sigman\}=1,\\  \exp\bigl(-\frac{2C}{\alpha_+(1-\sigman)}x^{1-\sigman}+\rO(x^{1-\min\{R,2\sigman\}})\bigr),& \text{if}\ \min\{R,2\sigman\}\neq1,\end{array}\right.$$
i.e., $\nu\in\cP^{2-}_{1<}$.

\noindent{$\bullet$} $R\ge\sigman=1$ and $\alpha=0$. If $R=\sigma$, then $$T_{\nu}(x)\gtrsim x^{-\frac{2C-\theta_{\nu}}{\alpha_+}},$$
i.e., $\nu\in\cP^{3-}$. If $R>\sigman$, then $$T_{\nu}(x)\gtrsim x^{-\frac{2C}{\alpha_+}},$$
which also indicates $\nu\in\cP^{3-}$.

\smallskip

\noindent{(vii)} $R=0$. From \eqref{Eq-5} it follows that \[1\ge\frac{T_{\nu}(x)}{T_{\nu}(x-1)}\ge\frac{\alpha_+-Cx^{-\sigman}}{\alpha_--\theta_{\nu}+Cx^{-\sigman}},\]
which yields $\frac{\alpha_+}{\alpha_--\theta_{\nu}}\le1$, i.e., $\alpha\le-\theta_{\nu}<0$. Similarly, based on \eqref{Eq-7-a}, $\theta_{\nu}\le\alpha_0\le\alpha_-$.

\smallskip

\noindent{$\bullet$} $R=0$, $\alpha+\theta_{\nu}=0$, $\sigman<1$:  $$T_{\nu}(x)\gtrsim\left\{\begin{array}{ll} \exp\bigl(-\frac{2C}{(\alpha_--\theta_{\nu})(1-\sigman)}x^{1-\sigman}+\rO(\log x)\bigr),& \text{if}\ 2\sigman=1,\\ \exp\bigl(-\frac{2C}{(\alpha_--\theta_{\nu})(1-\sigman)}x^{1-\sigman}+\rO(x^{1-2\sigman})\bigr),& \text{if}\ 2\sigman\neq1,\end{array}\right.$$
i.e., $\nu\in\cP^{2-}_{1-\sigman}$.

\noindent{$\bullet$} $R=0$, $\alpha+\theta_{\nu}=0$, $\sigman=1$:  $$T_{\nu}(x)\gtrsim x^{-\frac{2C}{\alpha_--\theta_{\nu}}},$$
i.e., $\nu\in\cP^{3-}$.

\noindent{$\bullet$} $R=0$, $\alpha+\theta_{\nu}<0$:  $$T_{\nu}(x)\gtrsim\left\{\begin{array}{ll} \exp\bigl(\frac{\alpha+\theta_{\nu}}{\alpha_--\theta_{\nu}} x+\rO(x^{1-\sigman})\bigr),& \text{if}\ \sigman<1,\\ \exp\bigl(\frac{\alpha+\theta_{\nu}}{\alpha_--\theta_{\nu}}x+\rO(\log x)\bigr),&\text{if}\ \sigman=1,\end{array}\right.$$
i.e., $\nu\in\cP^{2-}_{1}$.

\medskip

\noindent{\textbf{Step II}}. Establish upper estimates for $T_{\nu}$.

\smallskip

\noindent{Case I.} $R>1$. Similar arguments establishing the upper estimates for $T_{\pi}$ in the proof of Theorem~\ref{th-19} are adaptable to establish the those for $T_{\nu}$.

\medskip

\noindent{Latter part of (iv):} $R_->\max\{1,R_+\}$. Base on \eqref{Eq-8}, one can show there exists $C>0$ such that for all large $x$,

\begin{align*}
\frac{T_{\nu}(x)}{T_{\nu}(x-\omega_+)} &\le x^{R_+-R_-}\frac{\alpha_++Cx^{-\sigman}}{\alpha_{0}-\frac{\theta_{\nu}}{R_--1}(x-1)^{1-R_-}-Cx^{-\sigman}} \\
&=x^{R_+-R_-}
\lt(\tfrac{\alpha_+}{\alpha_{0}}+\rO(x^{-\min\{\sigman,R_--1\}})\rt),
\end{align*}
which implies that
\begin{align*}
T_{\nu}(x)&\lesssim \exp\Bigl(-(R_--R_+)\om_+^{-1}\log \Gamma(x\om_+^{-1})-\bigl((R_--R_+)\om_+^{-1}+\log\tfrac{\alpha_0}{\alpha_+}\bigr)x \\
&\quad +\rO(x^{1-\min\{\{\sigman,R_--1\}\}}+\log x)\Bigr).
\end{align*}
Then $\nu\in\cP^{1+}_{(R_--R_+)\om_+^{-1}}$.

\medskip

\noindent{\noindent{Latter part of (v)}:} $R_-=R_+>1$ and $\alpha<0$. An analogue of \eqref{SEq-43} is \begin{align}\nonumber &(\alpha_--\alpha_+)T_{\nu}(x)-\theta_{\nu}\sum_{y=x}^{\infty}y^{-R}{T_{\nu}(y)}+\sum_{j=\omega_-+1}^0\sum_{y=x}^{\infty}\nu(y)\beta_{j}(y+j)-
\sum_{j=1}^{\omega_+}\sum_{y=x}^{\infty}\nu(y)\beta_j(y+j)\\\label{SEq-43-b}
&=\sum_{j=\omega_-+1}^{-1}\sum_{\ell=0}^{-j-1}\nu\lt(x+\ell\rt)\widetilde{f}_{j}(x+j+\ell)+\sum_{j=1}^{\omega_+}
 \sum_{\ell=1}^{j}\nu\lt(x-\ell\rt)\widetilde{f}_{j}(x+j-\ell).
\end{align}
Based on this, one can show that there exists $C>0$ such that for all large $x$,
\begin{multline*}
$$\lt((\alpha_--\alpha_+)-\frac{\theta_{\nu}}{R-1}
x^{1-R}-\theta_{\nu}x^{-R}-Cx^{-\sigman}\rt)T_{\nu}(x)\le\text{LSH}_{\eqref{SEq-43-b}}\\ \le\lt(\alpha_--\alpha_+-\theta_{\nu}x^{-R}+Cx^{-\sigman}\rt)
T_{\nu}(x),\end{multline*}
\begin{multline*}(\alpha_+-Cx^{-\sigman})(T_{\nu}(x-1)-T_{\nu}(x))\le\text{RSH}_{\eqref{SEq-43-b}}\\ \le (\alpha_++Cx^{-\sigman})\lt(T_{\nu}(x-\omega_+)-T_{\nu}(x-\omega_--1)\rt).\end{multline*}
This implies
\[\frac{T_{\nu}(x-\om_+)}{T_{\nu}(x-\omega_--1)}\le\frac{\alpha_++Cx^{-\sigman}}{\alpha_--\frac{\theta_{\nu}}{R-1}
x^{1-R}-\theta_{\nu}x^{-R}},\]
and hence
\[T_{\nu}(x)\lesssim \exp\Bigl(-(\om_+-\om_--1)^{-1}\log\bigl(\frac{\alpha_-}{\alpha_+}\bigr)x+\rO(x^{1-\min\{\{\sigman,R-1\}\}}+\log x)\Bigr).\]
 Then $\nu\in\cP^{2+}_{1}$.

\smallskip

\noindent{Case II.} $R\le1$ (cases (viii)-(xi)).

Indeed, from \eqref{Eq-8}, for large $x$,
\begin{align*}\frac{T_{\nu}(x-\om_-)}{T_{\nu}(x)}&\le\frac{\alpha_-+Cx^{-\sigman}-\theta_{\nu}x^{-R}}{\alpha_-+Cx^{-\sigman}}=1-
\frac{\theta_{\nu}x^{-R}}{\alpha_-+Cx^{-\sigman}}\\
&=\left\{\begin{array}{ll} 1-\frac{\theta_{\nu}}{\alpha_-}x^{-R}+\rO(x^{-\min\{2R,R+\sigman\}}),& \text{if}\ R>0\\ 1-\frac{\theta_{\nu}}{\alpha_-}+\rO(x^{-\sigman}),& \text{if}\ R=0,\ \alpha_->\theta_{\nu},\\ \frac{C}{\alpha_-}x^{-\sigman}+\rO(x^{-2\sigman}), & \text{if}\ R=0,\ \alpha_-=\theta_{\nu}.\end{array}\right.\end{align*}
Using similar arguments as in the proof of Theorem~\ref{th-19}, we can show
\begin{align*}T_{\nu}(x)\lesssim\left\{\begin{array}{ll} x^{-\theta_{\nu}/\alpha_-},& \text{if}\ R=1,\\ \exp\lt(-\frac{\theta_{\nu}}{\alpha_-(1-R)}(-x\om_-^{-1})^{1-R}+\rO(\log x)\rt),& \text{if}\ 0<R<1,\\  &\quad\min\{2R,R+\sigman\}=1,\\ \exp\lt(-\frac{\theta_{\nu}}{\alpha_-(1-R)}(-x\om_-^{-1})^{1-R}+\rO(x^{1-\min\{2R,R+\sigman\}})\rt),& \text{if}\ 0<R<1,\\  &\quad \min\{2R,R+\sigman\}\neq1,\\
\exp\lt(\log\lt(1-\frac{\theta_{\nu}}{\alpha_-}\rt)(-x\om_-^{-1})+\rO(x^{1-\sigman})\rt),& \text{if}\ R=0,\ \alpha_->\theta_{\nu},\\
\Gamma(-x\om_-^{-1})^{-\sigman}(C\alpha_-^{-1})^{-x\om_-^{-1}-\widetilde{C}x^{1-\sigman}},& \text{if}\ R=0,\ \alpha_-=\theta_{\nu}.
\end{array}\right.\end{align*}
This implies that
\[\nu\in \left\{\begin{array}{ll} \cP^{3+}_{\theta_{\nu}/\alpha_-},& \text{if}\ R=1,\\ \cP^{2+}_{1-R},&\text{if}\ 0<R<1,\\ \cP^{2+}_{1},& \text{if}\ R=0,\ \alpha_->\theta_{\nu},\\
\cP^{1+}_{-\om_-^{-1}\sigman},& \text{if}\ R=0,\ \alpha_-=\theta_{\nu}.
\end{array}\right.\]



\bibliographystyle{plainnat}
{\footnotesize{\bibliography{references}}}

\appendix

\renewcommand\appendixname{Appendix}

\section{Structure of the state space}

\prob\label{pro-0}
Let $n=\min\cX$. Then $\cX\subseteq\om_*\N_0+n$. Assume $\rm{(\mathbf{A2})}$ and $\cT_+\neq\varnothing$. Then there exists $m\in\N_0$ with $m-n\in\om_*\N_0$ such that $\om_*\N_0+m\subseteq\partial^{\sf c}\subseteq\cX\subseteq\om_*\N_0+n$.
\proe
\prb
The first conclusion follows immediately from \cite[Lemma~B.2]{XHW20a}.

For the second conclusion, we only prove the case when $\cT_-=\varnothing$. The case when $\cT_-\neq\varnothing$ was proved in \cite[Lemma~B.1]{XHW20a}. If $\Omega_+=\{\om_*\}$ is a singleton, then by $\rm{(\mathbf{A2})}$, there exists $m\in\cX$ such that $\om_*\N_0+m\subseteq\cX$. Hence, the conclusion follows. Assume $\Omega_+$ has more than one element.
By the definition of $\om_*$, there exist coprime $\om_1,\ \om_2\in\cT_+$. By $\rm{(\mathbf{A2})}$, there exists $n_1\in\cX$ such that
$$q(x,x+\omega_1)>0,\ q(x,x+\omega_2)>0,\quad  x\in\cX,\ x\ge n_1.$$
Hence $(\om_1\N_0+n_1)\cup(\om_2\N_0+n_1)\subseteq\cX$. Further assume $n_1=0$ for the ease of exposition. We claim that there exists $s_j\in\cX\cap\N$ for $j=0,\ldots,\om_1-1$ such that $s_j-j\in\om_1\N_0$. Then \[\N_0+\prod_{j=0}^{\om_1-1}s_j\subseteq\cup_{j=0}^{\om_1-1}(\om_1\N_0+s_j)\subseteq\cX,\] since for every $x\in\N_0$, $x+\prod_{j=0}^{\om_1-1}s_j\in \om_1\N_0+(x\!\!\mod\om_1)$ and $\prod_{j=0}^{\om_1-1}s_j\ge s_k$ for all $k=0,\ldots,\om_1-1$. Hence it suffices to prove the above claim. Since $\om_1$ and $\om_2$ are coprime, there exist $m_1,\ m_2\in\N$ such that $m_1\om_1-m_2\om_2=1$. Then let $s_j=jm_1\om_1+m_2\om_2(\om_1-j)$, for $j=0,\ldots,\om_1-1$. It is ready to see that $s_j\in\om_1\N_0\cup\om_2\N_0\subseteq\cX$ and $$s_j=m_2\om_1\om_2+j(m_1\om_1-m_2\om_2)=m_2\om_2\om_1+j\in\om_1\N_0+j.$$ The proof is complete.
\pre

\section{A lemma for Theorem~\ref{th-19}}

\leb\label{Sle-1}
Assume $\rm{(\mathbf{A1})}$-$\rm{(\mathbf{A3})}$. Then
\begin{itemize}
\item[\rm(i)] $\alpha_0\ge\alpha_{-1}\ge \cdots\ge \alpha_{\omega_-+1}$, and $\alpha_1\ge \alpha_2\ge\cdots\ge \alpha_{\omega_+}$.
\item[\rm(ii)] $\beta_j(x)=\cab(\gamma_j+(-j+1)R_-\alpha_j)x^{-\sigman}+\rO(x^{-\rhon}),\quad \text{if}\ j=\om_-+1,\ldots,0,\\ (\gamma_j-jR_+\alpha_j)x^{-\sigman}+\rO(x^{-\rhon}),\qquad\qquad\, \text{if}\ j=1,\ldots,\om_+.\cae$
\item[\rm(iii)] It holds that   
\[\alpha_-=\omega_*\sum_{j=\omega_-+1}^0\alpha_{j},\ \alpha_+=\omega_*\sum_{j=1}^{\omega_+}\alpha_{j},\ \gamma_-=\omega_*\sum_{j=\omega_-+1}^0\gamma_{j},\ \gamma_+=\omega_*\sum_{j=1}^{\omega_+}\gamma_{j}.\]
\item[\rm(iv)]  If $\alpha=0$, then 
\[\sum_{j=\om_-+1}^{\omega_+}|j|\alpha_j=\vartheta\omega_*^{-2}.\]
\end{itemize}
\lee
\prb
Assume w.o.l.g. that $\omega_*=1$. We only prove the `+' cases. Analogous arguments apply to the `--' cases.

\smallskip
\noindent
{\rm(i)}-{\rm(ii)}. The first two properties follow directly from their definitions.

\smallskip
\noindent{\rm(iii)}. By Fubini's theorem,
\[\sum_{j=1}^{\omega_+}\sum_{\om\in A_j}q(x,x+\omega)=\sum_{j\ge1}\sum_{\ell\ge j}q(x,x+\ell)=\sum_{1\le j}jq(x,x+j)=\sum_{\omega\in\cT_+}q(x,x+\omega)\omega.\]
Comparing the coefficients before the highest degree of the polynomials on both sides, and use the definition of $\alpha_{\ell}$ as well as $\alpha_+$, we have
\[\alpha_+=\sum_{j=1}^{\om_+}\alpha_{j}.\]

\smallskip
\noindent
{\rm(iv)}. Due to Fubini's theorem again,
\[
\sum_{j=1}^{\omega_+}\sum_{\ell\ge j}jq(x,x+\ell)=\sum_{j=1}^{\omega_+}\frac{j(j+1)}{2}q(x,x+j).
\]
Note that $\alpha=0$ implies $\alpha_-=\alpha_+$. Since  $R_+=R$, comparing the coefficients before $x^R$, we have
\[\sum_{j=1}^{\omega_+}j\alpha_j=\frac{1}{2}\lt(\alpha_++\lim_{x\to\infty}\frac{\sum_{\omega\in\cT_+}
q(x,x+\omega)\omega^2}{x^{R_+}}\rt).\]
Similarly,
\[\sum_{j=\omega_-+1}^0\sum_{\ell\le j-1}|j-1|q(x,x+\ell)=\sum_{j=\om_-}^{-1}\frac{(j-1)j}{2}q(x,x+j).
\]
Hence
\[\sum_{j=\omega_-+1}^0|j-1|\alpha_j=\frac{1}{2}\lt(\alpha_-+\lim_{x\to\infty}\frac{\sum_{\omega\in\cT_-}
q(x,x+\omega)\omega^2}{x^{R_-}}\rt).\]
Then the conclusion follows from \[\vartheta=\frac{1}{2}\lim_{x\to\infty}\frac{\sum_{\omega\in\cT}
q(x,x+\omega)\omega^2}{x^{R}}\]
and
\[\sum_{j=\om_-+1}^{0}|j-1|\alpha_j=\alpha_-+\sum_{j=\om_-+1}^{0}|j|\alpha_j.\]
\pre

\section{Calculation of sharp asymptotics of stationary distributions in special cases}
We first provide the sharp asymptotics of stationary distributions for BDPs. For any two real-valued functions $f,g$ on $\R$, we write $f(x)\sim g(x)$ if there exists $C>1$ such that $C^{-1}g(x)\le f(x)\le Cg(x)$ for all $x\in\R$.

\prob\label{pro-15}
Assume ($\rm\mathbf{A1}$)-($\rm\mathbf{A3}$), $\alpha=0$, and $\om_+=-\om_-=1$.
 Let $\pi$ be a stationary distribution of $Y_t$ supported on $\cX$. Then

\begin{itemize}
\item[(i)] if $\Delta=0$, $\sigman<1$, and $\min\{2\sigman,\rhon\}>1$, then $R>1$ and $\pi\in\cP^{3+}_{R-1}\cap \cP^{3-}_{R-1}$.
\item[(ii)] if $\Delta>0$, $\sigman<1$,
then $\pi\in\cP^{3+}_{1-\sigman}\cap \cP^{3-}_{1-\sigman}$.
\item[(iii)] if $\sigman=1$, then $\Delta>1$ and $\pi\in\cP^{3+}_{\Delta-1}\cap \cP^{3-}_{\Delta-1}$.
\end{itemize}
\proe

{\begin{proof}
The proof is similar to that of Theorem~\ref{th-19}, with the aid of Stirling's formula and Euler-Maclaurin's formula.\end{proof}}

When $Y_t$ is not a BDP, the asymptotic tail of a stationary distribution can be established in some cases when $\alpha=0$. When $\alpha=0$, $\alpha_+=\alpha_-$ and $\frac{1}{\om_+-\om_--1}\le\frac{2}{\om_+-\om_-}\le\frac{\alpha_+\om_*}{\vartheta}\le1$. Hence $$\delta\le R-\frac{\gamma}{\vartheta}=\Delta\frac{\alpha_+\om_*}{\vartheta}\le\Delta,$$ and both equalities hold if and only if $Y_t$ is a BDP.

\prob\label{pro-16}
Assume ($\rm\mathbf{A1}$)-($\rm\mathbf{A3}$), $\alpha=0$, $\gamma+\vartheta<0$, and $\partial=\varnothing$. Let $\pi$ be a stationary distribution of $Y_t$ supported on $\cX$.
Then for large $x\in\cX$,
\[ T_{\pi}(x)\sim \cab  \exp(\frac{\gamma}{\vartheta(1-\sigman)}(\om_*^{-1}x)^{1-\sigman}+\rO(\log x)),\ \hspace{2.8cm} \text{if}\ \sigman<1,\\ \hspace{8.85cm}  \min\{2\sigman,\rhon\}=1,\\
x^{1-R+\sigman}\exp(\frac{\gamma}{\vartheta(1-\sigman)}(\om_*^{-1}x)^{1-\sigman}+\rO(x^{1-\min\{2\sigman,\rhon\}})),\ \text{if}\ \sigman<1,\\ \hspace{8.85cm} \min\{2\sigman,\rhon\}\neq1,\\
x^{1+\frac{\gamma}{\vartheta}-R}\lt(1+\rO\lt(x^{-1}\rt)\rt),\  \hspace{4.7cm}\text{if}\ \sigman=1.\cae\]
 Hence $\pi\in\cP^{2-}_{1-\sigman}\cap\cP^{2+}_{1-\sigman}$ if $\sigman<1$, and $\pi\in\cP^{3-}_{R-\frac{\gamma}{\vartheta}-1}\cap\cP^{3+}_{R-\frac{\gamma}{\vartheta}-1}$ if $\sigman=1$.
\proe
{\begin{proof}
 To prove the conclusions, we   apply \cite[Theorem\,1]{DKW13} and  arguments similar to those in the proof of Theorem~\ref{th-19}.
\end{proof}}

\end{document}